\documentstyle{amsppt}

\nologo

\refstyle{A}

\font\bigbold=cmbx12

\define\imdf{\operatorname{im}d\tilde f}

\define\Fdot{\bold F^\bullet}

\define\Pdot{\bold P^\bullet}

\define\con{\overline{T^*_{{}_{S_\alpha}}\Cal U}}

\define\conc{\Big[ \overline{T^*_{{}_{S_\alpha}}\Cal U}\Big ]}

\define\dm{\operatorname{dim}}

\define\xreg{X_{\operatorname{reg}}}

\define\kp{{}^k\hskip -.01in\Pdot}

\hoffset .1 true in \voffset .2 true in

\hsize=6.3 true in \vsize=8.7 true in

\topmatter

\title Critical Points of Functions on Singular Spaces \endtitle

\author David B. Massey \endauthor

\address{David B. Massey, Dept. of Mathematics, Northeastern University,  Boston, MA, 02115, USA}
\endaddress

\email{DMASSEY\@neu.edu}\endemail



\abstract We compare and contrast various notions of the ``critical locus'' of a complex analytic function on a singular
space. After choosing a topological variant as our primary notion of the critical locus, we justify our choice by
generalizing L\^e and Saito's result that constant Milnor number implies that Thom's $a_f$ condition is satisfied.
\endabstract

\endtopmatter

\document

\vskip .3in

\noindent\S0. {\bf Introduction}.  

\vskip .2in

Let $\Cal U$ be an open subset of $\Bbb C^{n+1}$, let $\bold z:= (z_0, z_1, \dots, z_n)$ be coordinates for $\Bbb
C^{n+1}$, and suppose that
$\tilde f:\Cal U\rightarrow\Bbb C$ is an analytic function. Then, all conceivable definitions of the {\it critical
locus},
$\Sigma \tilde f$, of $\tilde f$ agree: one can consider the points, $\bold x$, where the derivative vanishes, i.e.,
$d_\bold x \tilde f=0$, or one can consider the points,
$\bold x$, where the Taylor series of $\tilde f$ at $\bold x$ has no linear term, i.e., $\tilde f-\tilde f(\bold
x)\in\frak m^2_{{}_{\Cal U,
\bold x}}$ (where $\frak m_{{}_{\Cal U, \bold x}}$ is the maximal ideal in the coordinate ring of $\Cal U$ at
$\bold x$), or one can consider the points, $\bold x$, where the Milnor fibre of $\tilde f$ at $\bold x$,
$F_{\tilde f,\bold x}$, is not trivial (where, here, ``trivial'' could mean even up to analytic isomorphism). 

Now, suppose that $X$ is an analytic subset of $\Cal U$, and let $f:= \tilde f_{|_{X}}$. Then, what should be meant by
``the critical locus of $f$''? It is not clear what the relationship is between points, $\bold x$, where $f-f(\bold
x)\in\frak m^2_{{}_{X,
\bold x}}$ and points where the Milnor fibre, $F_{f, \bold x}$, is not trivial (with any definition of trivial);
moreover, the derivative
$d_\bold x f$ does not even exist.

\vskip .1in

We are guided by the successes of Morse Theory and stratified Morse Theory to choosing the Milnor fibre definition as
our primary notion of critical locus, for we believe that critical points should coincide with changes in the topology
of the level hypersurfaces of $f$. Therefore, we make the following definition: 

\vskip .3in

\noindent{\bf Definition 0.1}. The {\it $\Bbb C$-critical locus of $f$}, $\Sigma_{{}_\Bbb C}f$, is given by 
$$\Sigma_{{}_\Bbb C}f:= \{\bold x\in X\ |\  H^*(F_{f, \bold x};\ \Bbb C)\neq H^*(point;\ \Bbb C)\}.$$
 (The reasons for using field coefficients, rather than $\Bbb Z$, are
technical: we want Lemma 3.1 to be true.)

\vskip .3in

In Section 1, we will compare and contrast the $\Bbb C$-critical locus with other possible notions of critical locus,
including the ones mentioned above and the stratified critical locus.

\vskip .1in

After Section 1, the remainder of this paper is dedicated to showing that Definition 0.1 really yields a useful,
calculable definition of the critical locus. We show this by looking at the case of a generalized isolated singularity,
i.e., an isolated point of $\Sigma_{{}_{\Bbb C}}f$, and showing that, at such a point, there is a workable definition
of the Milnor number(s) of $f$; we show that the Betti numbers of the Milnor fibre can be calculated (3.7.ii), and we
give a generalization of the result of L\^e and Saito [{\bf L-S}] that constant Milnor number throughout a family
implies Thom's $a_f$ condition holds. Specifically, in Corollary 5.14, we prove (with slightly weaker hypotheses) that:

\vskip .3in

\noindent{\bf Theorem 0.2}. {\it   Let $W$ be a (not necessarily purely) $d$-dimensional analytic subset of an open subset of $\Bbb C^n$. Let
$Z$ be a $d$-dimensional irreducible component of $W$. Let $X:=\overset\circ\to{\Bbb D}\times W$ be the product of an open disk about the
origin with $W$, and let $Y:=\overset\circ\to{\Bbb D}\times Z$.

Let $f:(X, \overset\circ\to{\Bbb D}\times\{\bold 0\})\rightarrow (\Bbb C, 0)$ be an analytic function, and let $f_t(\bold z):= f(t, \bold
z)$. Suppose that $f_0$ is in the square of the maximal ideal of $Z$ at $\bold 0$.

Suppose that $\bold 0$ is an isolated point of $\Sigma_{{}_{\Bbb
C}}(f_0)$, and that the reduced Betti number $\tilde b_{d-1}(F_{f_a, (a,\bold 0)})$ is independent of $a$ for all
small $a$. 

Then, $\tilde b_{d-1}(F_{f_a, (a,\bold 0)})\neq 0$ and, near
$\bold 0$, 
$\Sigma(f_{|_{Y_{\operatorname{reg}}}})\subseteq \overset\circ\to{\Bbb D}\times\{\bold 0\}$ and
the pair $(Y_{\operatorname{reg}} -\Sigma (f_{|_{Y_{\operatorname{reg}}}}),\ \overset\circ\to{\Bbb D}\times\{\bold 0\})$ satisfies Thom's
$a_f$ condition at $\bold 0$. 

}

\vskip .3in

Thom's $a_f$ is important for several reasons, but perhaps the best reason is because it is an hypothesis of Thom's
Second Isotopy Lemma. General results on the $a_f$ condition have been proved by many researchers: Hironaka, L\^e,
Saito, Henry, Merle, Sabbah, Brian\c con, Maisonobe, Parusi\'nski, etc., and the above theorem is closely related to
the recent results contained in [{\bf BMM}] and [{\bf P}]. However, the reader should contrast the hypotheses of Theorem 0.2 with
those of the main theorem of [{\bf BMM}] (Theorem 4.2.1); our main hypothesis is that a single number is constant throughout the
family, while the main hypothesis of Theorem 4.2.1 of [{\bf BMM}] is a condition which requires one to check an infinite amount of
data: the property of local stratified triviality. Moreover, the Betti numbers that we require to be constant are actually
calculable.

\vskip .3in

While much of this part is fairly technical in nature, there are three new, key ideas that guide us throughout.

\vskip .1in

The first of these fundamental precepts is: controlling the {\bf vanishing cycles} in a family of functions is enough to
control Thom's $a_f$ condition and, perhaps, the topology throughout the family. While this may seem like an obvious
principle -- given the results of L\^e and Saito in [{\bf L-S}] and of L\^e and Ramanujam in [{\bf L-R}] -- in fact, in
the general setting, most of the known results seem to require the constancy of much stronger data, e.g., the constancy
of the polar multiplicities [{\bf Te}] or that one has the local stratified triviality property [{\bf BMM}]. In a very
precise sense, controlling the polar multiplicities corresponds to controlling the {\bf nearby cycles} of the family of
functions, instead of merely controlling the vanishing cycles. As we show in Corollary 4.4, controlling the
characteristic cycle of the vanishing cycles is sufficient for obtaining the $a_f$ condition. 

Our second fundamental idea is: the correct setting for all of our cohomological results is where perverse sheaves are
used as coefficients. While papers on intersection cohomology abound, and while perverse sheaves are occasionally used
as a tool (e.g., [{\bf BMM}, 4.2.1]), we are not aware of any other work on general singularities in which arbitrary
perverse sheaves of coefficients are used in an integral fashion throughout. The importance of perverse
sheaves in this paper begins with Theorem 3.2, where we give a description of the critical locus of a function with
respect to a perverse sheaf.

The third new feature of this paper is the recurrent use of the perverse cohomology of a complex of sheaves. This device
allows us to take our general results about perverse sheaves and translate them into statements about the constant
sheaf.  The reason that we use perverse cohomology, instead of intersection cohomology, is because perverse cohomology
has such nice functorial properties: it commutes with Verdier dualizing, and with taking nearby and vanishing cycles
(shifted by $[-1]$). If we were only interested in proving results for local complete intersections (l.c.i.'s), we would
never need the perverse cohomology; however, we want to prove completely general results. The perverse cohomology seems
to be a hitherto unused tool for accomplishing this goal.

\vskip .3in

\noindent This paper is organized as follows:

\vskip .1in

In Section 1, we discuss seven different notions of the ``critical locus'' of a function. We give examples
to show that, in general, all of these notions are different.

Section 2 is devoted to proving an ``index theorem'', Theorem 2.10, which provides the main link between the
topological data of the Milnor fibre and the algebraic data obtained by blowing-up the image of $d\tilde f$ inside the
appropriate space. This theorem is presented with coefficients in a bounded, constructible complex of sheaves; this
level of generality is absolutely necessary in order to obtain the results in the remainder of this paper.

Section 3 uses the index theorem of Section 2 to show that $\Sigma_{{}_\Bbb C}f$ and the Betti numbers of the Milnor
fibre really are fairly well-behaved. This is accomplished by applying Theorem 2.10 in the case where the complex of
sheaves is taken to be the perverse cohomology of the shifted constant sheaf. Perverse cohomology essentially gives us
the ``closest'' perverse sheaf to the constant sheaf. Many of the results of Section 3 are stated for arbitrary perverse
sheaves, for this seems to be the most natural setting.

Section 4 contains the necessary results from conormal geometry that we will need in order to conclude that topological
data implies that Thom's $a_f$ condition holds. The primary result of this section is Corollary 4.4, which once again
relies on the index theorem from Section 2.

Section 5 begins with a discussion of ``continuous families of constructible complexes of sheaves''. We then prove in Theorem 5.7 that
additivity of Milnor numbers occurs in  continuous families of perverse sheaves, and we use this to conclude additivity of the Betti numbers
of the Milnor fibres, by once again resorting to the perverse cohomology of the shifted constant sheaf. Finally, in Corollaries 5.11
and 5.12, we prove that the constancy of the Milnor/Betti number(s) throughout a family implies that the $a_f$ condition
holds -- we prove this first in the setting of arbitrary perverse sheaves, and then for perverse cohomology of the
shifted constant sheaf. By translating our hypotheses from the language of the derived category back into more
down-to-Earth terms, we obtain Corollary 5.12, which leads to Theorem 0.2 above.

\vskip .5in

\noindent{\bigbold Section 1}. {\bf CRITICAL AVATARS}.  

\vskip .2in

We continue with $\Cal U$, $\bold z$, $\tilde f$, $X$, and $f$ as in the introduction.

\vskip .1in

In this section, we will investigate seven possible notions of the ``critical locus'' of a function on a
singular space, one of which is the $\Bbb C$-critical locus already defined in 0.1.

\vskip .3in

\noindent{\bf Definition 1.1}. The {\it algebraic critical locus of $f$}, $\Sigma_{\operatorname{alg}}f$, is defined by
$$
\Sigma_{\operatorname{alg}}f := \{\bold x\in X\ |\ f-f(\bold x)\in\frak m^2_{X, \bold x}\} .
$$ 

\vskip .3in

\noindent{\it Remark 1.2}. It is a trivial exercise to verify that 
$$
\Sigma_{\operatorname{alg}}f  = \{\bold x\in X\ |\ \text{there exists a local extension, $\hat f$, of $f$ to
$\Cal U$ such that $d_\bold x\hat f = 0$}\}.
$$

Note that $\bold x$ being in $\Sigma_{\operatorname{alg}}f$ does {\bf not} imply that {\bf every} local extension of $f$
has zero for its derivative at $\bold x$. 

One might expect that $\Sigma_{\operatorname{alg}}f$ is always a closed set; in fact, it need not be. Consider the
example where $X:=V(xy)\subseteq\Bbb C^2$, and $f = y_{|_X}$. We leave it as an exercise for the reader to verify that
$\Sigma_{\operatorname{alg}}f = V(y)-\{\bold 0\}$.

\vskip .3in

There are five more variants of the critical locus of $f$ that we will consider. We let $\xreg$ denote the regular (or
smooth) part of $X$ and, if $M$ is an analytic submanifold of $\Cal U$, we let
$T^*_{{}_M}\Cal U$ denote the conormal space to $M$ in $\Cal U$ (that is, the elements $(\bold x, \eta)$ of the
cotangent space to $\Cal U$ such that $\bold x\in M$ and $\eta$ annihilates the tangent space to $M$ at $\bold x$). We
let $N(X)$ denote the Nash modification of $X$, so that the fibre $N_{\bold x}(X)$ at $\bold x$ consists of limits of
tangent planes from the regular part of $X$.

We
also remind the reader that complex analytic spaces possess canonical Whitney stratifications (see [{\bf Te}]).

\vskip .3in

\noindent{\bf Definition 1.3}. We define the {\it regular critical locus of $f$}, $\Sigma_{\operatorname{reg}}f$, to be
the critical locus of the restriction of $f$ to $\xreg$, i.e., $\Sigma_{\operatorname{reg}}f =
\Sigma\big(f_{|_{\xreg}}\big)$.

\vskip .3in

We define the {\it Nash critical locus of $f$}, $\Sigma_{{}_{\operatorname{Nash}}}f$, to be 
$$
\big\{\bold x\in X\ |\ \text{there exists a local extension, $\hat f$, of $f$ to $\Cal
U$ such that $d_{\bold x}\hat f(T)\equiv 0$, for all $T\in N_{\bold x}(X)$}\big\}.
$$

\vskip .3in

We define the {\it conormal-regular critical locus of $f$}, $\Sigma_{\operatorname{cnr}}f$, to be 
$$
\big\{\bold x\in X\ |\ \text{there exists a local extension, $\hat f$, of $f$ to $\Cal
U$ such that $(\bold x, d_\bold x\hat f)\in \overline{T^*_{{}_{\xreg}}\Cal U}$ }\big\};
$$
it is trivial to see that this set is equal to
$$
\big\{\bold x\in X\ |\ \text{there exists a local extension, $\hat f$, of $f$ to $\Cal
U$ such that $d_{\bold x}\hat f(T)\equiv 0$, for some $T\in N_{\bold x}(X)$}\big\}.
$$

\vskip .3in

Let $\Cal S= \{S_\alpha\}$ be a (complex analytic) Whitney stratification of $X$. We define the {\it $\Cal
S$-stratified critical locus of $f$}, $\Sigma_{\Cal S}f$, to be $\bigcup_\alpha\Sigma\big(f_{|_{S_\alpha}}\big)$. If $\Cal S$ is
clear, we simply call $\Sigma_{\Cal S}f$  the stratified critical locus.

If $\Cal S$ is, in fact, the canonical Whitney stratification of $X$, then we write $\Sigma_{\operatorname{can}}f$ in place
of $\Sigma_{\Cal S}f$, and call it the {\it  canonical stratified critical locus}.

\vskip .3in

We define the {\it relative differential critical locus of $f$}, $\Sigma_{\operatorname{rdf}}f$, to be the union of the singular
set of $X$ and $\Sigma_{\operatorname{reg}}f$. 

If $\bold x\in X$ and $h_1, \dots, h_j$ are equations whose zero-locus defines $X$
near $\bold x$, then $\bold x\in \Sigma_{\operatorname{rdf}}f$ if and only if the rank of the Jacobian map of $(\tilde f, h_1,
\dots, h_j)$ at $\bold x$ is not maximal among all points of $X$ near $\bold x$. By using this Jacobian, we could (but will
not) endow
$\Sigma_{\operatorname{rdf}}f$ with a scheme structure (the {\it critical space}) which is independent of the choice of the
extension
$\tilde f$ and the defining functions $h_1, \dots, h_n$ (see [{\bf Lo}, 4.A]). The proof of the independence uses relative
differentials; this is the reason for our terminology.

\vskip .3in

\noindent{\it Remark 1.4}. In terms of conormal geometry, $\Sigma_{\Cal S}f = \Big\{\bold x\in X\ |\ (\bold x, d_\bold x\tilde f)\in
\bigcup_\alpha T^*_{S_\alpha}\Cal U\Big\}$ or, using Whitney's condition a) again, $\Sigma_{\Cal S}f = \Big\{\bold x\in
X\ |\ (\bold x, d_\bold x\tilde f)\in \bigcup_\alpha\con\Big\}$.

\vskip .1in

Clearly, $\Sigma_{\operatorname{rdf}}f$ is closed, and it is an easy exercise to show that Whitney's
condition a) implies that $\Sigma_{\Cal S}f$ is closed. On the other hand, $\Sigma_{\operatorname{reg}}f$ is, in general, not
closed and, in order to have any information at singular points of $X$, we will normally look at its closure
$\overline{\Sigma_{\operatorname{reg}}f}$.

Looking at the definition of $\Sigma_{\operatorname{cnr}}f$, one might expect that
$\overline{\Sigma_{\operatorname{reg}}f}=\Sigma_{\operatorname{cnr}}f$. In fact, we shall see in Example 1.8 that this
is false. That $\Sigma_{\operatorname{cnr}}f$ is, itself, closed is part of the following proposition. (Recall that
$\tilde f$ is our fixed extension of $f$ to all of $\Cal U$.)

\vskip .4in

In the following proposition, we show that, in the definitions of the Nash and conormal-regular critical loci, we could have used
``for all'' in place of ``there exists'' for the local extensions; in particular, this implies that we can use the fixed extension
$\tilde f$. Finally, we show that the conormal-regular critical locus is closed.

\vskip .3in

\noindent{\bf Proposition 1.5}. {\it The Nash critical locus of $f$ is equal to
$$
\big\{\bold x\in X\ |\ \text{for all local extensions, $\hat f$, of $f$ to $\Cal
U$, $d_{\bold x}\hat f(T)\equiv 0$, for all $T\in N_{\bold x}(X)$}\big\}=$$
$$ \big\{\bold x\in X\ |\ d_{\bold x}\tilde f(T)\equiv 0, \text{ for all }T\in N_{\bold x}(X)\big\}.
$$

\vskip .2in

The conormal-regular critical locus of $f$ is equal to
$$
\big\{\bold x\in X\ |\ \text{for all local extensions, $\hat f$, of $f$ to $\Cal U$,
$(\bold x, d_\bold x\hat f)\in \overline{T^*_{{}_{\xreg}}\Cal U}$ }\big\}$$
$$ =\big\{\bold x\in X\ |\ (\bold x, d_\bold x\tilde f)\in \overline{T^*_{{}_{\xreg}}\Cal U} \big\}.
$$

In addition, $\Sigma_{\operatorname{cnr}}f$ is closed.

}

\vskip .2in

\noindent{\it Proof}. Let $Z:=\big\{\bold x\in X\ |\ \text{for all local extensions, $\hat f$, of $f$ to $\Cal
U$, $d_{\bold x}\hat f(T)\equiv 0$, for all $T\in N_{\bold x}(X)$}\big\}$. Clearly, we have $Z\subseteq
\Sigma_{{}_{\operatorname{Nash}}}f$. 

Suppose now that $\bold x\in\Sigma_{{}_{\operatorname{Nash}}}f$. Then, there exists a local extension, $\hat f$, of $f$ to
$\Cal U$ such that $d_{\bold x}\hat f(T)\equiv 0$, for all $T\in N_{\bold x}(X)$.  Let $\check f$ be another local extension of $f$
to $\Cal U$ and let $T_\infty\in N_{\bold x}(X)$; to show that $\bold x\in Z$, what we must show is that  $d_{\bold x}\check
f(T_\infty)\equiv 0$. 

Suppose not.  Then, there exists $\bold v\in T_\infty$ such that $d_\bold x\check f(\bold v)\neq 0$, but $d_\bold x\hat f(\bold
v)= 0$. Therefore, there exist
$\bold x_i\in X_{\operatorname{reg}}$ and $\bold v_i\in T_{\bold x_i}X_{\operatorname{reg}}$ such that $\bold x_i\rightarrow \bold x$,
$T_{\bold x_i}X_{\operatorname{reg}}\rightarrow T_\infty$, and $\bold v_i\rightarrow \bold v$.

Let $\Cal V$ be an open neighborhood of $x$ in $\Cal U$ which in $\hat f$ and $\check f$ are both defined. Let $\Phi:\Cal
V\cap\overline{TX_{\operatorname{reg}}}\rightarrow\Bbb C$ be defined by $\Phi(\bold p, \bold w)= d_\bold p(\hat f-\check f)(\bold
w)$. Then, $\Phi$ is continuous, and so $\Phi^{-1}(0)$ is closed. As $(\hat f-\check f)_{|_X}\equiv 0$, $(\bold x_i, \bold v_i)\in
\Phi^{-1}(0)$, and thus $(\bold x, \bold v)\in \Phi^{-1}(0)$ -- a contradiction. Therefore,
$Z=\Sigma_{{}_{\operatorname{Nash}}}f$.

\vskip .2in

It follows immediately that $\Sigma_{{}_{\operatorname{Nash}}}f = \big\{\bold x\in X\ |\ d_{\bold x}\tilde f(T)\equiv 0, \text{ for
all }T\in N_{\bold x}(X)\big\}$. 

\vskip .3in

Now, let $W:=\big\{\bold x\in X\ |\ \text{for all local extensions, $\hat f$, of $f$ to $\Cal U$,
$(\bold x, d_\bold x\hat f)\in \overline{T^*_{{}_{\xreg}}\Cal U}$ }\big\}$. Clearly, we have $W\subseteq
\Sigma_{\operatorname{cnr}}f$. 

Suppose now that $\bold x\in\Sigma_{\operatorname{cnr}}f$. Then, there exists a local extension, $\hat f$, of $f$ to
$\Cal U$ such that $(\bold x, d_\bold x\hat f)\in \overline{T^*_{{}_{\xreg}}\Cal U}$. Let $(\bold x_i, \eta_i)\in
T^*_{{}_{\xreg}}\Cal U$ be such that $(\bold x_i,
\eta_i)\rightarrow (\bold x, d_\bold x\hat f)$. Let $\check f$ be another local extension of $f$ to $\Cal U$; to show
that $\bold x\in W$, what we must show is that  $(\bold x, d_\bold x\check f)\in
\overline{T^*_{{}_{\xreg}}\Cal U}$.

Since $(\check f-\hat f)_{|_X}\equiv 0$, for all $\bold q\in \xreg$, $\big(\bold q, d_\bold q(\check f-\hat f)\big)\in
T^*_{{}_{\xreg}}\Cal U$; in particular, $\big(\bold x_i, d_{\bold x_i}(\check f-\hat f)\big)\in T^*_{{}_{\xreg}}\Cal U$.
Thus,  $\big(\bold x_i, \eta_i+d_{\bold x_i}(\check f-\hat f)\big)\in T^*_{{}_{\xreg}}\Cal U$, and $\big(\bold x_i,
\eta_i+d_{\bold x_i}(\check f-\hat f)\big)\rightarrow (\bold x, d_\bold x\check f)$.  Therefore, $(\bold x, d_\bold
x\check f)\in
\overline{T^*_{{}_{\xreg}}\Cal U}$, and $W=\Sigma_{\operatorname{cnr}}f$. 

\vskip .2in

It follows immediately that $\Sigma_{\operatorname{cnr}}f = \big\{\bold x\in X\ |\ (\bold x, d_\bold x\tilde f)\in
\overline{T^*_{{}_{\xreg}}\Cal U} \big\}$.  

\vskip .2in

Finally, we need to show that $\Sigma_{\operatorname{cnr}}f$ is closed. Let $\Psi:X\rightarrow T^*\Cal U$ be given by
$\Psi(\bold x)=(\bold x, d_{\bold x}\tilde f)$. Then, $\Psi$ is a continuous map and, by the above,
$\Sigma_{\operatorname{cnr}}f=\Psi^{-1}(\overline{T^*_{{}_{\xreg}}\Cal U})$.
\qed

\vskip .3in

\noindent{\bf Proposition 1.6}. {\it There are inclusions
$$
\overline{\Sigma_{\operatorname{reg}}f}\ \subseteq\
\overline{\Sigma_{\operatorname{alg}}f}\ \subseteq\ \overline{\Sigma_{{}_{\operatorname{Nash}}}f}\ \subseteq\
\Sigma_{\operatorname{cnr}}f\
\subseteq\
\overline{\Sigma_{{}_\Bbb C}f}\ \subseteq\ \Sigma_{\operatorname{can}}f\ \subseteq\ \Sigma_{\operatorname{rdf}}f.
$$
In addition, if $\Cal S$ is a Whitney stratification of $X$, then $\Sigma_{\operatorname{can}}f\ \subseteq\ \Sigma_{\Cal S}f$
}.

\vskip .3in

\noindent{\it Proof}. Clearly, $\Sigma_{\operatorname{reg}}f\ \subseteq\
\Sigma_{\operatorname{alg}}f\ \subseteq\
\Sigma_{{}_{\operatorname{Nash}}}f\ \subseteq\
\Sigma_{\operatorname{cnr}}f$, and so the containments for their closures follows (recall, also that $\Sigma_{\operatorname{cnr}}f$
is closed). It is also obvious that
$\Sigma_{\operatorname{can}}f\ \subseteq\ \Sigma_{\operatorname{rdf}}f$ and 
$\Sigma_{\operatorname{can}}f\ \subseteq\ \Sigma_{\Cal S}f$.

That $\Sigma_{{}_\Bbb C}f\ \subseteq\ \Sigma_{\operatorname{can}}f$ follows from Stratified Morse Theory [{\bf
Go-Mac1}], and so, since $\Sigma_{\operatorname{can}}f$ is closed, $\overline{\Sigma_{{}_\Bbb C}f}\ \subseteq\
\Sigma_{\operatorname{can}}f$.

It remains for us to show that $\Sigma_{\operatorname{cnr}}f\ \subseteq\
\overline{\Sigma_{{}_\Bbb C}f}$. Unfortunately, to reach this conclusion, we must refer ahead to Theorem 3.6, from which
it follows immediately. (However, that $\overline{\Sigma_{\operatorname{alg}}f}\ \subseteq\
\overline{\Sigma_{{}_\Bbb C}f}$ follows from A'Campo's Theorem [{\bf A'C}].)\qed

\vskip .3in

\noindent{\it Remark 1.7}. For a fixed stratification $\Cal S$, for all $\bold x\in X$, there exists a neighborhood
$\Cal W$ of $\bold x$ in $X$ such that $\Cal W\cap\Sigma_{\Cal S}f\subseteq V(f-f(\bold x))$. This is easy to show: the
level hypersurfaces of $f$ close to $V(f-f(\bold x))$ will be transverse to all of the strata of $\Cal S$ near $\bold
x$. All of our other critical loci which are contained in $\Sigma_{\Cal S}f$ (i.e., all of them except
$\Sigma_{\operatorname{rdf}}f$)  also satisfy this local isolated critical value property.

\vskip .3in

\noindent{\it Example 1.8}. In this example, we wish to look at the containments given in Proposition 1.6, and
investigate whether the containments are proper, and also investigate what would happen if we did not take closures in
the four cases where we do.

\vskip .1in

The same example that we used in Remark 1.2 shows that none of $\Sigma_{\operatorname{reg}}f$,
$\Sigma_{\operatorname{alg}}f$, $\Sigma_{{}_{\operatorname{Nash}}}f$, or $\Sigma_{{}_\Bbb C}f$ are necessarily closed; if
$X:=V(xy)\subseteq\Bbb C^2$, and $f = y_{|_X}$, then all four critical sets are precisely $V(y)-\{\bold 0\}$.
Additionly, since
$\Sigma_{\operatorname{cnr}}f = V(y)$, this example also shows that, in general, 
$\Sigma_{\operatorname{cnr}}f\ \not\subseteq\ \Sigma_{{}_\Bbb C}f$.

\vskip .1in

If we continue with $X=V(xy)$ and let $g:=(x+y)^2_{|_X}$, then $\overline{\Sigma_{\operatorname{alg}}g}=\{\bold 0\}$ and
$\overline{\Sigma_{\operatorname{reg}}g}=\emptyset$; thus, in general, $\overline{\Sigma_{\operatorname{reg}}f}\neq
\overline{\Sigma_{\operatorname{alg}}f}$.

\vskip .1in

While it is easy to produce examples where $\Sigma_{{}_{\operatorname{Nash}}}f$  is not equal to $\Sigma_{\operatorname{alg}}f$ and
examples where $\Sigma_{{}_{\operatorname{Nash}}}f$  is not equal to $\Sigma_{\operatorname{cnr}}f$, it is not quite so easy to come
up with examples where all three of these sets are distinct. We give such an example here.

Let
$Z:=V((y-zx)(y^2-x^3))\subseteq\Bbb C^3$ and
$L:= y_{|_Z}$. Then, one easily verifies that $\Sigma_{\operatorname{alg}}f=\emptyset$, $\Sigma_{{}_{\operatorname{Nash}}}f=\{\bold
0\}$, and $\Sigma_{\operatorname{cnr}}f=\Bbb C\times\{\bold 0\}$.

\vskip .1in

If  $X=V(xy)$ and  $h:=(x+y)_{|_X}$, then $\Sigma_{{}_\Bbb C}h = \{\bold 0\}$ and
$\Sigma_{\operatorname{cnr}}h=\emptyset$; thus, in general,  $\Sigma_{\operatorname{cnr}}f\neq\overline{\Sigma_{{}_\Bbb
C}f}$.

\vskip .1in

Let $W:= V(z^5+ty^6z+y^7x+x^{15})\subseteq \Bbb C^4$; this is the example of Brian\c con and Speder [{\bf B-S}] in which
the topology along the $t$-axis is constant, despite the fact that the origin is a point-stratum in the canonical
Whitney stratification of $W$. Hence, if we let $r$ denote the restriction of $t$ to $W$, then, for values of $r$ close
to $0$, $\bold 0$ is the only point in $\Sigma_{\operatorname{can}}r$ and $\bold 0\not\in \Sigma_{{}_\Bbb C}r$.
Therefore, $\bold 0\in
\Sigma_{\operatorname{can}}r\ -\ \overline{\Sigma_{{}_\Bbb C}r}$, and so, in general, $\overline{\Sigma_{{}_\Bbb
C}f}\neq\Sigma_{\operatorname{can}}f$.

\vskip .1in

Using the coordinates $(x,y,z)$ on $\Bbb C^3$, consider the cross-product $Y:=V(y^2-x^3)\subseteq\Bbb C^3$. The canonical Whitney
stratification of $Y$ is given by $\{Y-\{\bold 0\}\times\Bbb C, \ \{\bold 0\}\times\Bbb C\}$. Let $\pi:=z_{|_{Y}}$. Then,
$\Sigma_{\operatorname{can}}\pi=\emptyset$, while $\Sigma_{\operatorname{rdf}}\pi = \{\bold 0\}\times\Bbb C$. Thus, in general, 
$\Sigma_{\operatorname{can}}f\neq\Sigma_{\operatorname{rdf}}f$.

\vskip .1in

It is, of course, easy to throw extra, non-canonical, Whitney strata into almost any example in order to see that, in
general, 
$\Sigma_{\operatorname{can}}f\neq\Sigma_{\Cal S}f$.

\vskip .3in

To summarize the contents of this example and Proposition 1.6: we have seven seemingly reasonable definitions of
``critical locus'' for complex analytic functions on singular spaces (we are not counting $\Sigma_{\Cal S}f$, since it
is not intrinsically defined). All of our critical locus avatars agree for manifolds. The sets
$\Sigma_{\operatorname{reg}}f$,
$\Sigma_{\operatorname{alg}}f$, $\Sigma_{{}_{\operatorname{Nash}}}f$, and $\Sigma_{{}_\Bbb C}f$ need not be closed. There is a chain
of containments among the closures of these critical loci, but -- in general -- none of the sets are equal.

\vskip .3in

However, we consider the sets 
$\overline{\Sigma_{\operatorname{reg}}f}$, $\overline{\Sigma_{\operatorname{alg}}f}$, $\overline{\Sigma_{{}_{\operatorname{Nash}}}f}$,
and
$\Sigma_{\operatorname{cnr}}f$ to be too small; these ``critical loci'' do not detect the change in topology at the level
hypersurface $h=0$ in the simple example $X=V(xy)$ and 
$h=(x+y)_{|_X}$ (from Example 1.8).

Despite the fact that the Stratified Morse Theory of [{\bf Go-Mac1}] yields nice results and requires one to
consider the stratified critical locus, we also will not use $\Sigma_{\operatorname{can}}f$ (or any other $\Sigma_{\Cal
S}f$) as our primary notion of critical locus;
$\Sigma_{\operatorname{can}}f$ is often too big. As we saw in the Brian\c con-Speder example in Example 1.8, the
stratified critical locus sometimes forces one to consider ``critical points'' which do not correspond to changes in
topology.

Certainly, $\Sigma_{\operatorname{rdf}}f$ is far too large, if we want critical points to have {\bf any} relation to changes in
the topology of level hypersurfaces: if $X$ has a singular set $\Sigma X$, then the critical space of the
projection $\pi: X\times \Bbb C\rightarrow\Bbb C$ would consist of $\Sigma X\times\Bbb C$, despite the obvious triviality of the
family of level hypersurfaces defined by $\pi$.

\vskip .1in

Therefore, we choose to concentrate our attention on the $\Bbb C$-critical locus, and we will justify this choice with the
results in the remainder of this paper.

 Note that we consider $\Sigma_{{}_\Bbb C}f$, not its closure, to be the correct
notion of critical locus; we think that this is the more natural definition, and we consider the question of when
$\Sigma_{{}_\Bbb C}f$ is closed to be an interesting one. It is true, however, that all of our results refer to
$\overline{\Sigma_{{}_\Bbb C}f}$. We should mention here that, while $\Sigma_{{}_\Bbb C}f$ need not be closed, the
existence of Thom stratifications [{\bf Hi}] implies that $\Sigma_{{}_\Bbb C}f$ is at least analytically constructible;
hence,
$\overline{\Sigma_{{}_\Bbb C}f}$ is an analytic subset of $X$.

\vskip .3in

Before we leave this section, in which we have already looked at seven definitions of ``critical locus'', we
need to look at one last variant. As we mentioned at the end of the introduction, even though we wish to investigate the
Milnor fibre with coefficients in $\Bbb C$, the fact that the shifted constant sheaf on a non-l.c.i. need not be
perverse requires us to take the perverse cohomology of the constant sheaf.  This means that we need to consider the
hypercohomology of Milnor fibres with coefficients in an arbitrary bounded, constructible complex of sheaves (of $\Bbb
C$-vector spaces).

The $\Bbb C-$critical locus is  nicely described in terms of vanishing cycles (see [{\bf K-S}] for general properties of
vanishing cycles, but be aware that we use the more traditional shift):
$$
\Sigma_{{}_{\Bbb C}}f = \{\bold x\in X\ |\ H^*(\phi_{f-f(\bold x)}\Bbb C^\bullet_X)_\bold x\neq 0\}.
$$  This definition generalizes easily to yield a definition of the critical loci of $f$ with respect to arbitrary
bounded, constructible complexes of sheaves on $X$.

\vskip .2in

Let $\Cal S:=\{S_\alpha\}$ be a Whitney stratification of $X$, and let $\Fdot$ be a bounded complex of sheaves (of $\Bbb
C$-vector spaces) which is constructible with respect to $\Cal S$.

\vskip .3in

\noindent{\bf Definition 1.9}. The {\it $\Fdot$-critical locus  of $f$} , $\Sigma_{{}_{\Fdot}}f$, is defined by
$$
\Sigma_{{}_{\Fdot}}f := \{\bold x\in X\ |\ H^*(\phi_{f-f(\bold x)}\Fdot)_\bold x\neq 0\}.
$$ 

\vskip .2in

\noindent{\it Remark 1.10}. Stratified Morse Theory (see [{\bf Go-Mac1}]) implies that $\Sigma_{{}_{\Fdot}}f
\ \subseteq\ \Sigma_{\Cal S}f$ (alternatively, this follows from 8.4.1 and 8.6.12 of [{\bf K-S}], combined with the
facts that complex analytic Whitney stratifications are $w$-stratifications, and  $w$-stratifications are
$\mu$-stratifications.)

\vskip .3in

We could discuss three more notions of the critical locus of a function -- two of which are obtained by
picking specific complexes for $\Fdot$ in Definition 1.9. However, we will defer the introduction of these new critical
loci until Section 3; at that point, we will have developed the tools necessary to say something interesting about
these three new definitions.

\vskip .5in

\noindent{\bigbold Section 2}. {\bf THE LINK BETWEEN THE ALGEBRAIC AND\newline\hbox{}\hskip .93in TOPOLOGICAL POINTS OF
VIEW}.

\vskip .2in

We continue with our previous notation:  $X$ is a complex analytic space contained in some open subset $\Cal U$ of some
$\Bbb C^{n+1}$, $\tilde f:\Cal U\rightarrow\Bbb C$ is a complex analytic function, $f=\tilde f_{|_X}$, $\Cal
S=\{S_\alpha\}$ is a Whitney stratification of $X$ with connected strata, and  $\Fdot$ is a bounded complex of sheaves
(of $\Bbb C$-vector spaces) which is constructible with respect to $\Cal S$. In addition, $N_\alpha$ and $\Bbb
L_\alpha$ are, respectively, the normal slice and complex link of the $d_\alpha$-dimensional stratum $S_\alpha$ (see [{\bf Go-Mac1}]).

	In this section, we are going to prove a general result  which describes the characteristic cycle of $\phi_f\Fdot$ in
terms of blowing-up the image of $d\tilde f$ inside the conormal spaces to strata. We will have to wait until the next
section (on results for perverse sheaves) to actually show how this provides a relationship between
$\Sigma_{{}_{\Fdot}}f$ and 
$\Sigma_{\Cal S}f$ in the case where $\Fdot$ is perverse.

Beginning in this section, we will use some aspects of intersection theory, as described in [{\bf F}]; however, at all
times, the setting for our intersections will be the most trivial: we will only consider proper intersections of complex
analytic cycles (not cycle classes) inside an ambient analytic manifold. In this setting, there is a well-defined
intersection cycle.

\vskip .2in

\noindent{\bf Definition 2.1}. Recall that the {\it characteristic cycle,  $\operatorname{Ch}(\Fdot)$, of $\Fdot$} in
$T^*\Cal U$ is the linear combination 
$\sum_{\alpha} m_{\alpha}(\Fdot) \conc$, where   the $m_\alpha(\Fdot)$ are integers given by  
$$m_\alpha(\Fdot)\ :=  \ (-1)^{\dm X-1}\chi(\phi_{L_{|_X}}\bold F^\bullet)_{\bold x} \ = \  (-1)^{\dm
X-d_\alpha-1}\chi(\phi_{L_{|_{N_\alpha}}}{\bold F^\bullet}_{|_{N_\alpha}})_{\bold x}$$  for any point $\bold x$ in 
$S_\alpha$, with normal slice $N_\alpha$  at $\bold x$,  and any  $L: (\Cal U, x) \rightarrow\ (\Bbb C,0)$
such that  $d_{\bold x}L$  is a non-degenerate covector at $\bold x$ (with respect to our fixed stratification; see
[{\bf Go-Mac1}])  and
$L_{|_{S_\alpha}}$ has a Morse singularity at $\bold x$.  This cycle is independent of all the choices made (see, for
instance, [{\bf K-S}, Chapter IX]).

\vskip .4in

We need a number of preliminary results before we can prove the main theorem (Theorem 2.10) of this section.

\vskip .2in

\noindent{\bf Definition 2.2}. Recall that, if $M$ is an analytic submanifold of $\Cal U$ and $M\subseteq X$, then the
{\it  relative conormal space (of $M$ with respect to $f$ in $\Cal U$)}, $T^*_{f_{|_M}}\Cal U$, is given by
$$ T^*_{f_{|_M}}\Cal U:=\{(\bold x, \eta)\in T^*\Cal U\ |\ \bold x\in M,\ \eta\big(\operatorname{ker}d_\bold
x(f_{|_M})\big)=0\} = 
$$
$$
\{(\bold x, \eta)\in T^*\Cal U\ |\ \bold x\in M,\ \eta\big(T_\bold x M\cap\operatorname{ker}d_\bold x\tilde f\big)=0\}.
$$

We define the {\it total relative conormal cycle}, $T^*_{{}_{f, \Fdot}}\Cal U$, by $\dsize T^*_{{}_{f, \Fdot}}\Cal U :=
\sum_{S_\alpha\not\subseteq f^{-1}(0)}m_\alpha \Big [\overline{T^*_{f_{|_{S_\alpha}}}\Cal U}\Big ]$.
\vskip .3in

\vskip .5in

{\bf From this point, through Lemma 2.9, it will be convenient to assume that we have refined our stratification
$\Cal S = \{S_\alpha\}$  so that $V(f)$ is a union of strata. By Remark 1.7, this implies that, in a neighborhood of
$V(f)$, if $S_\alpha\not\subseteq V(f)$, then $\Sigma(f_{|_{S_\alpha}})=\emptyset$. }

\vskip .5in

We shall need the following important result from [{\bf BMM}, 3.4.2].  

\vskip .3in

\noindent{\bf Theorem 2.3}. ([{\bf BMM}]) {\it The characteristic cycle of the sheaf of nearby cycles of $\Fdot$ along
$f$,
$\operatorname{Ch}\big(\psi_f\Fdot\big)$, is isomorphic to the intersection product
$T^*_{{}_{f, \Fdot}}\Cal U\cdot \big(V(f)\times
\Bbb C^{n+1}\big)$ in $\Cal U\times \Bbb C^{n+1}$. }

\vskip .3in

 Let $\Gamma^1_{f, L}(S_\alpha)$ denote the closure in
$X$ of the relative polar curve of $f$ with respect to $L$ inside $S_\alpha$ (see [{\bf M1}] and [{\bf M3}]).  It is
important to note that $\Gamma^1_{f, L}(S_\alpha)$ is the closure of the polar curve in $S_\alpha$, {\bf not} in 
$\overline{S_\alpha}$; that is, $\Gamma^1_{f, L}(S_\alpha)$ has no components contained in any strata
$S_\beta\subseteq\overline{S_\alpha}$ such that $S_\beta\neq S_\alpha$.

It is convenient to have a specific point in $X$ at which to work. Below, we concentrate our attention at the origin; of
course, if the origin is not in $X$ (or, if the origin is not in $V(f)$), then we obtain zeroes for all the terms below.
For any bounded, constructible complex
$\bold A^\bullet$ on a subspace of $\Cal U$, let $m_{\bold 0}(\bold A^\bullet)$ equal the coefficient of
$\Big[T^*_{\{\bold 0\}}\Cal U\Big]$ in the characteristic cycle of $\bold A^\bullet$. 

We need to state one further result without proof  -- this result can be obtained from [{\bf BMM}], but we give the
result as stated in [{\bf M1}, 4.6].

\vskip .3in

\noindent{\bf Theorem 2.4}.  {\it  For  generic linear forms $L$,  we have the following formulas:

\vskip .1in

$$ m_\bold 0(\psi_f\Fdot)=\sum_{S_\alpha\not\subseteq V(f)}m_\alpha 
\big(\Gamma^1_{{}_{f, L}}(S_\alpha)\cdot V(f)\big)_\bold 0 ;
$$

\vskip .1in

$$ m_\bold 0(\Fdot)+m_\bold 0(\Fdot_{|_{V(f)}})=\sum_{S_\alpha\not\subseteq V(f)}m_\alpha 
\big(\Gamma^1_{{}_{f, L}}(S_\alpha)\cdot V(L)\big)_\bold 0 ; \text{ and}
$$ 

\vskip .1in

$$ m_\bold 0(\phi_f\Fdot)=m_\bold 0(\Fdot) + \sum_{S_\alpha\not\subseteq V(f)}m_\alpha 
\left(\big(\Gamma^1_{{}_{f, L}}(S_\alpha)\cdot V(f)\big)_\bold 0 - 
\big(\Gamma^1_{{}_{f, L}}(S_\alpha)\cdot V(L)\big)_\bold 0\right) . 
$$ }

\vskip .4in

\noindent{\bf Lemma 2.5}. {\it If $S_\alpha\not\subseteq f^{-1}(0)$, then the coefficient of
$\left[\Bbb P(T^*_{\{\bold 0\}}\Cal U)\right]$ in $\Bbb P\big(\overline{T^*_{{}_{f_{|_{S_\alpha}}}}\Cal U}\big)\cdot
\big(V(f)\times \Bbb P^n\big)$ is given by 
$\big(\Gamma^1_{{}_{f, L}}(S_\alpha)\cdot V(f)\big)_\bold 0$. }

\vskip .3in

\noindent{\it Proof}.  Take a complex of sheaves, $\bold F^\bullet$, which has a characteristic cycle consisting only of
$\conc$  (see, for instance, [{\bf M1}]).  Now, apply the formula for
$m_\bold 0(\psi_f\Fdot)$ from Theorem 2.4 together with Theorem 2.3. \qed

\vskip .4in

We need to establish some notation that we shall use throughout the remainder of this section.

\vskip .2in

 Using the isomorphism,
$T^*\Cal U\cong
\Cal U\times\Bbb C^{n+1}$, we consider
$\operatorname{Ch}(\Fdot)$ as a cycle in
$X\times \Bbb C^{n+1}$; we use $\bold z:= (z_0, \dots, z_n)$ as coordinates on $\Cal U$ and $\bold w:= (w_0, \dots,
w_n)$ as the cotangent coordinates.  

Let $I$ denote the sheaf of ideals on $\Cal U$ given by the image of $d\tilde f$, i.e., 
$I=\big<w_0-\frac{\partial \tilde f}{\partial z_0}, \dots, w_n-\frac{\partial
\tilde f}{\partial z_n}\big>$.  For all $\alpha$, let $B_\alpha =
\operatorname{Bl}_{\operatorname{im} d\tilde f}\con$ denote the blow-up of
$\con$ along the image of $I$ in $\con$, and let $E_\alpha$ denote the corresponding exceptional divisor.  For all
$\alpha$, we have
$E_\alpha\subseteq B_\alpha\subseteq X\times\Bbb C^{n+1}\times \Bbb P^n$.  Let $\pi: X\times\Bbb C^{n+1}\times \Bbb
P^n\rightarrow X\times \Bbb P^n$ denote the projection.  Note that, if $(\bold x, \bold w, [\eta])\in E_\alpha$, then
$\bold w = d_\bold x\tilde f$ and so, for all
$\alpha$,
$\pi$ induces an isomorphism from $E_\alpha$ to
$\pi(E_\alpha)$.  We refer to $E := \sum_\alpha m_\alpha E_\alpha$ as the {\it total exceptional divisor} inside the
{total blow-up}
$\operatorname{Bl}_{\operatorname{im} d\tilde f}\operatorname{Ch}(\Fdot) := 
\sum_\alpha m_\alpha\operatorname{Bl}_{\operatorname{im} d\tilde f}\conc$.

\vskip .4in

\noindent{\bf Lemma 2.6}. {\it For all $S_\alpha$, there is an inclusion
$\pi\left(\operatorname{Bl}_{\imdf}\con\right)\subseteq 
\Bbb P\big(\overline{T^*_{{}_{f_{|_{S_\alpha}}}}\Cal U}\big)$. }

\vskip .3in

\noindent{\it Proof}.  This is entirely straightforward.  Suppose that
$$(\bold x, \bold w, [\eta])\in 
\operatorname{Bl}_{\imdf}\con =
\overline{\operatorname{Bl}_{\imdf}T^*_{{}_{S_\alpha}}\Cal U} .$$  Then, we have a sequence $(\bold
x_i, \bold w_i, [\eta_i])\in 
\operatorname{Bl}_{\imdf}T^*_{{}_{S_\alpha}}\Cal U$ such that $(\bold x_i,
\bold w_i, [\eta_i])\rightarrow (\bold x, \bold w, [\eta])$.  

By definition of the blow-up, for each $(\bold x_i, \bold w_i, [\eta_i])$, there exists a sequence
$(\bold x^j_i,
\bold w^j_i)\in T^*_{{}_{S_\alpha}}\Cal U - \imdf$ such that $(\bold x^j_i,
\bold w^j_i, [\bold w^j_i - d_{\bold x^j_i}\tilde f])\rightarrow (\bold x_i,
\bold w_i, [\eta_i])$.  Now, $(\bold x^j_i, [\bold w^j_i - d_{\bold x^j_i}\tilde f])$ is clearly in
$\Bbb P\big(T^*_{{}_{f_{|_{S_\alpha}}}}\Cal U\big)$, and so each $(\bold x_i, [\eta_i])$ is in
$\Bbb P\big(\overline{T^*_{{}_{f_{|_{S_\alpha}}}}\Cal U}\big)$.  Therefore, 
$(\bold x, [\eta])\in \Bbb P\big(\overline{T^*_{{}_{f_{|_{S_\alpha}}}}\Cal U}\big)$. \qed

\vskip .4in

\noindent{\bf Lemma 2.7}. {\it If $S_\alpha\not\subseteq f^{-1}(0)$, then the coefficient of
$\big[\Bbb P\big(T^*_{{}_{\{\bold 0\}}}\Cal U\big)\big]=\{\bold 0\}\times\Bbb P^n$ in $\pi_*(E_\alpha)$ equals 
$\big(\Gamma^1_{{}_{f, L}}(S_\alpha)\cdot V(f)\big)_\bold 0 - 
\big(\Gamma^1_{{}_{f, L}}(S_\alpha)\cdot V(L)\big)_\bold 0$.
 }

\vskip .3in

\noindent{\it Proof}. We will work inside $\Cal U\times\Bbb C^{n+1}\times\Bbb P^n$. We use $[u_0:\dots : u_n]$ as projective
coordinates, and calculate the coefficient of 
$G_{\bold 0}:=\big[\{\bold 0\}\times \{d_{\bold 0}\tilde f\}\times\Bbb P^n\big]$ in $E_\alpha$ using the affine patch
$\{u_0 \neq 0\}$.  

Letting $\tilde u_i = u_i/u_0$ for
$i\geqslant 1$, we have 
$\{u_0 \neq 0\} \cap B_\alpha=$ 
$$\overline{\left\{\big(\bold x, \bold w, (\tilde u_1, \dots, \tilde u_n)\big)\in \big(\con -
\operatorname{im} d\tilde f\big)\times \Bbb C^n \
\Big| \ w_i - \frac{\partial\tilde f}{\partial z_i} = \tilde u_i\left(w_0 -
\frac{\partial\tilde f}{\partial z_0}\right)\ , i\geqslant 1\right\}} ,$$ and $\{u_0 \neq 0\} \cap E_\alpha$ equals the
intersection product $\big(\{u_0
\neq 0\}
\cap B_\alpha\big)\cdot V\left(w_0 - \frac{\partial\tilde f}{\partial z_0}\right)$ in $\Cal U\times\Bbb
C^{n+1}\times\Bbb C^n$.

\vskip .1in

To calculate the multiplicity of $\{u_0 \neq 0\} \cap G_\bold 0$ in  
$\Big(\{u_0 \neq 0\} \cap B_\alpha\Big) \ \cdot \ V\left(w_0 - \frac{\partial\tilde f}{\partial z_0}\right)$, we move to
a generic point of
$\{u_0 \neq 0\} \cap
G_\bold 0$ and take a normal slice; that is, we fix a generic choice 
$(\tilde u_1, \dots, \tilde u_n) = (a_1, \dots, a_n)$.  This corresponds to choosing the generic linear form $L = z_0 +
a_1z_1+\dots a_nz_n$.

\vskip .3in

We claim that $Z:=\{u_0\neq 0\}\cap B_\alpha\cap V(\tilde u_1-a_1, \dots, \tilde u_n-a_n) - \{\bold 0\}\times\Bbb
C^{n+1}\times\Bbb C^n$  equals the set of all $\big(\bold x, \bold w, (a_1, \dots, a_n)\big)$ such that $\bold x\in
\Gamma^1_{{}_{f, L}}(S_\alpha) - \{\bold 0\}$ and $\bold w = d_\bold x\tilde f- \lambda(\bold x)d_\bold x L$, where
$\lambda(\bold x)$ is the unique non-zero complex number such that $\big(d_\bold x\tilde f- \lambda(\bold x)d_\bold x
L\big)(T_\bold x S_\alpha) = 0$.  

Once we show that $\bold x$ must be in $\Gamma^1_{{}_{f, L}}(S_\alpha) - \{\bold 0\}$, then it follows at once from the
definition of the relative polar curve that there exists a $\lambda(\bold x)$ as above.  That such a $\lambda(\bold x)$
must be unique is easy:  if we had two distinct such
$\lambda$, then we would have 
$d_\bold x L(T_\bold x S_\alpha) = 0$ -- but this is impossible for generic $L$.

Now, by definition of the relative polar curve and using 2.6, we find 
$$
\pi\big(\{u_0\neq 0\}\cap B_\alpha\cap V(\tilde u_1-a_1, \dots, \tilde u_n-a_n) -
\{\bold 0\}\times\Bbb C^{n+1}\times\Bbb C^n\big)
$$
$$
\subseteq\{u_0\neq 0\}\cap \Bbb P\big(\overline{T^*_{{}_{f_{|_{S_\alpha}}}}\Cal U}\big) 
\cap V(\tilde u_1-a_1, \dots, \tilde u_n-a_n) -
\{\bold 0\}\times\Bbb C^{n+1}\times\Bbb C^n
$$
$$ = \big(\Gamma^1_{{}_{f, L}}(S_\alpha)-\{\bold 0\}\big)\times\{(a_1, \dots, a_n)\} .
$$ (Actually, here we have also used 2.4 and 2.5 to conclude that there are no components of the relative polar curve
which are contained in strata other than $S_\alpha$.)

Thus, 
$$Z = \big\{(\bold x, \bold w, (a_1, \dots, a_n)) \ \big | \ \bold x\in
 \big(\Gamma^1_{{}_{f, L}}(S_\alpha)-\{\bold 0\}\big) \text{ and } 
\bold w = d_\bold x\tilde f- \lambda(\bold x)d_\bold x L\big\} ,$$ and the coefficient of
$G_\bold 0$ in $E_\alpha$ equals the intersection
number $$\left(\overline{Z}\cdot V\Big(w_0 - \frac{\partial\tilde f}{\partial z_0}\Big)\right)_{(\bold 0, \bold 0, (a_1,
\dots, a_n))}$$ in $\Cal U\times\Bbb C^{n+1}\times\{(a_1, \dots, a_n)\}$.

\vskip .2in

Now, for each component $C$ of $\overline{Z}$ through $(\bold 0, \bold 0, (a_1, \dots, a_n))$, select a local analytic
parameterization $\bold u_{{}_C}(t) = (\bold x_{{}_C}(t),
\bold w_{{}_C}(t), (a_1,
\dots, a_n))\in C$ such that $\bold x_{{}_C}(0) = \bold 0$, $\bold w_{{}_C}(0) = \bold 0$, and, for
$t\neq 0$,
$\bold u_{{}_C}(t)\in C-\{(\bold 0, \bold 0)\}\times \Bbb C^n$.  Then, 
$$\left(\overline{Z}\cdot V\Big(w_0 - \frac{\partial\tilde f}{\partial z_0}\Big)\right)_{(\bold 0,
\bold 0, (a_1, \dots, a_n))} = \ \sum_C \operatorname{mult}\Big\{
\Big( w_0 - \frac{\partial\tilde f}{\partial z_0}\Big)\circ \bold u_{{}_C}(t)\Big\} .$$ Moreover, a quick look at the
definition of $Z$ tells us that 
$\Big( w_0 - \frac{\partial\tilde f}{\partial z_0}\Big)\circ \bold u_{{}_C}(t) = \lambda(\bold x_{{}_C}(t))$. Thus, what
we want to show is that 
$$
\sum_C \operatorname{mult}\lambda(\bold x_{{}_C}(t)) = 
\big(\Gamma^1_{{}_{f, L}}(S_\alpha)\cdot V(f)\big)_\bold 0 - 
\big(\Gamma^1_{{}_{f, L}}(S_\alpha)\cdot V(L)\big)_\bold 0 .
$$

If we look now at $\big(\Gamma^1_{{}_{f, L}}(S_\alpha)\cdot V(f)\big)_\bold 0$, we find 
$$\big(\Gamma^1_{{}_{f, L}}(S_\alpha)\cdot V(f)\big)_\bold 0= 
\sum_C \operatorname{mult}f(\bold x_{{}_C}(t)) =
 \sum_C\Big(1 + \operatorname{mult}\big(f(\bold x_{{}_C}(t))\big)^\prime\Big)$$
$$ = \sum_C\Big(1 + \operatorname{mult}d_{{}_{\bold x_{{}_C}(t)}}\tilde f(\bold x^\prime_{{}_C}(t))
 \Big) \  = \ \sum_C\Big(1 + \operatorname{mult}\big((\bold w_{{}_C}(t) + \lambda(\bold x_{{}_C}(t)) d_{{}_{\bold
x_{{}_C}(t)}}L)\circ(\bold x^\prime_{{}_C}(t))
 \big)\Big).
$$ As $(\bold x_{{}_C}(t), \bold w_{{}_C}(t))\in T^*_{{}_{S_\alpha}}\Cal U$ for $t\neq 0$,
$\bold w_{{}_C}(t)\circ \bold x^\prime_{{}_C}(t) = 0$.   In addition,  
$d_{{}_{\bold x_{{}_C}(t)}}L\circ\bold x^\prime_{{}_C}(t) = \big(L(\bold x_{{}_C}(t))\big)^\prime$, and so we obtain
that 
$$\big(\Gamma^1_{{}_{f, L}}(S_\alpha)\cdot V(f)\big)_\bold 0 =$$
$$\sum_C
\big(\operatorname{mult}L(p_{{}_C}(t)) + \operatorname{mult}\lambda(p_{{}_C}(t))\big) = 
\big(\Gamma^1_{{}_{f, L}}(S_\alpha)\cdot V(L)\big)_\bold 0 +
\sum_C \operatorname{mult}\lambda(\bold x_{{}_C}(t)) ,$$ and so we are finished.  \qed

\vskip .4in

\noindent{\bf Lemma 2.8}. {\it For all $\alpha$ such that $S_\alpha\subseteq V(f)$, there is an inclusion of the
exceptional divisor 

$$E_\alpha\cong \pi(E_\alpha)\subseteq \Bbb P\big(\overline{T^*_{{}_{f_{|_{S_\alpha}}}}\Cal U}\big)\cap
\big(V(f)\times \Bbb P^n\big) .$$ }

\vskip .3in

\noindent{\it Proof}.  That $\pi$ is an isomorphism when restricted to the exceptional divisor is trivial:  $(\bold x,
\bold w, [\eta])\in E_\alpha$ implies that $\bold w = d_\bold x\tilde f$.  From Lemma 2.6, $\pi(E_\alpha)\subseteq 
\pi\left(\operatorname{Bl}_{\imdf}\con\right)\subseteq 
\Bbb P\big(\overline{T^*_{{}_{f_{|_{S_\alpha}}}}\Cal U}\big)$. The result follows.
\qed

\vskip .4in

\noindent{\bf Lemma 2.9}. {\it If $S_\alpha\subseteq f^{-1}(0)$, then
$E_\alpha\cong\pi(E_\alpha) = \Bbb P(\con)$.
 }

\vskip .3in

\noindent{\it Proof}.  If $S_\alpha\subseteq f^{-1}(0)$, then  
$\Bbb P\big(T^*_{{}_{f_{|_{S_\alpha}}}}\Cal U\big) = \Bbb P\big(T^*_{{}_{S_\alpha}}\Cal U\big)$, and so, by 2.8, 
$\pi(E_\alpha)\subseteq \Bbb P(\con)$.  We will demonstrate the reverse inclusion.

Suppose that we have $(\bold x, [\eta])\in \Bbb P(\con)$.  Then, there exists a sequence 
$(\bold x_i, \eta_i)\in T^*_{{}_{S_\alpha}}\Cal U$ such that $(\bold x_i,
\eta_i)\rightarrow (\bold x, \eta)$. Hence, $\big(\bold x_i, \ 
\frac{1}{i}\eta_i+d_{\bold x_i}\tilde f\big)\in
\con-\imdf$ and 
$$\Big(\bold x_i, \ 
\frac{1}{i}\eta_i+d_{\bold x_i}\tilde f, \
\Big[\Big(\frac{1}{i}\eta_i+d_{\bold x_i}\tilde f\Big) - d_{\bold x_i}\tilde f\Big]\Big)\rightarrow (\bold x, d_\bold x
\tilde f, [\eta])\in E_\alpha . \qed$$

\vskip .5in

We come now to the main theorem of this section. This theorem relates the topological data provided by the vanishing cycles of a
function $f$ to the algebraic data given by blowing-up the image of the differential of an extension of $f$.

\vskip .3in

\noindent{\bf Theorem 2.10}. {\it The projection $\pi$ induces an isomorphism between the total exceptional divisor  
$E\subseteq\operatorname{Bl}_{\operatorname{im} d\tilde f}\operatorname{Ch}(\Fdot)$ and the sum over all $v\in\Bbb C$ of
the  projectivized characteristic cycles of the sheaves of vanishing cycles of $\Fdot$ along $f-v$, i.e., 
$$E\cong \pi_*(E) = \sum_{v\in\Bbb C}\Bbb P(\operatorname{Ch}(\phi_{f-v}\Fdot)).$$
 } 

\vskip .2in

\noindent{\it Proof}.  Remarks 1.7 and 1.10 imply that, locally, $\operatorname{supp}\phi_{f-v}\Fdot\subseteq
f^{-1}(v)$. As the
$\Bbb P(\operatorname{Ch}(\phi_{f-v}\Fdot))$ are disjoint for different values of $v$, we may immediately reduce
ourselves to the case where we are working near $\bold 0\in X$ and where $f(\bold 0)=0$. We refine our stratification so
that, for all
$\alpha$,
$\Sigma(f_{|_{S_\alpha}})=\emptyset$ unless $S_\alpha\subseteq V(f)$.  As any newly introduced stratum will appear with
a coefficient of zero in the characteristic cycle, the total exceptional divisor will not change. We need to show that
$E\cong \pi(E) = \Bbb P(\operatorname{Ch}(\phi_{f}\Fdot))$.

\vskip .3in

\noindent Now, we will first show that $\pi(E)$ is Lagrangian.  

\vskip .2in

If $S_\alpha\subseteq f^{-1}(0)$, then
$\pi(E_\alpha) = \Bbb P(\con)$ by 2.9. If $S_\alpha\not\subseteq f^{-1}(0)$, then, by Theorem 2.3, $\Bbb
P\big(\overline{T^*_{{}_{f_{|_{S_\alpha}}}}\Cal U}\big)\cap
\big(V(f)\times \Bbb P^n\big)$ is Lagrangian and, in particular, is purely
$n$-dimensional.  By Lemma 2.8, $\pi(E_\alpha)$ is a  purely $n$-dimensional analytic set contained in  $\Bbb
P\big(\overline{T^*_{{}_{f_{|_{S_\alpha}}}}\Cal U}\big)\cap
\big(V(f)\times \Bbb P^n\big)$.  We need to show that $\pi(E_\alpha)$ is closed.

Suppose we have a sequence $(\bold x_i, [\eta_i])\in \pi(E_\alpha)$ and $(\bold x_i, [\eta_i])\rightarrow (\bold x,
[\eta])$ in $\Cal U\times\Bbb P^n$.  Then, there exists a sequence
$\bold w_i$ so that $(\bold x_i, \bold w_i, [\eta_i])\in E_\alpha$; by definition of the exceptional divisor, this
implies $\bold w_i = d_{\bold x_i}\tilde f$.  Therefore, $(\bold x_i,
\bold w_i, [\eta_i])\rightarrow (\bold x, d_{\bold x}\tilde f, [\eta])$, which is contained in
$E_\alpha$ since $E_\alpha$ is closed in $\Cal U\times\Bbb C^{n+1}\times\Bbb P^n$.  Thus, 
$(\bold x, [\eta])\in \pi(E_\alpha)$, and so 
$\pi(E_\alpha)$ is closed and, hence, Lagrangian.

\vskip .2in

Now, $\pi(E)$ and $\Bbb P(\operatorname{Ch}(\phi_f\Fdot))$ are both supported over $\Sigma_{\Cal S}f$ and, by taking
normal slices to strata, we are reduced to the point-stratum case.  Thus, what we need to show is: the coefficient of
$\big[\Bbb P\big(T^*_{{}_{\{\bold 0\}}}\Cal U\big)\big]$ in $E$ equals the coefficient of
$\big[\Bbb P\big(T^*_{{}_{\{\bold 0\}}}\Cal U\big)\big]$ in 
$\Bbb P(\operatorname{Ch}(\phi_f\Fdot))$.  Using 2.4, this is equivalent to showing that the coefficient of 
$\big[\Bbb P\big(T^*_{{}_{\{\bold 0\}}}\Cal U\big)\big]$ in $E$ equals 
$$m_\bold 0(\Fdot) + \sum_{S_\alpha\not\subseteq V(f)}m_\alpha 
\left(\big(\Gamma^1_{{}_{f, L}}(S_\alpha)\cdot V(f)\big)_\bold 0 - 
\big(\Gamma^1_{{}_{f, L}}(S_\alpha)\cdot V(L)\big)_\bold 0\right)$$ for a generic linear form $L$.

But, by 2.9, 
$$E \ = \ \sum_\alpha m_\alpha E_\alpha \ = \ \sum_{{}_{S_\alpha\subseteq V(f)}}m_\alpha \big[\Bbb P(\con)\big] \ + \ 
\sum_{{}_{S_\alpha\not\subseteq V(f)}}m_\alpha E_\alpha$$ and the coefficient of $\big[\Bbb P\big(T^*_{{}_{\{\bold
0\}}}\Cal U\big)\big]$ in $\dsize\sum_{{}_{S_\alpha\subseteq V(f)}}m_\alpha { \big[\Bbb P(\con)\big]}$ is precisely
$m_\bold 0(\Fdot)$.

\vskip .3in

Therefore, we will be finished if we can show that the coefficient of 
$\big[\Bbb P\big(T^*_{{}_{\{\bold 0\}}}\Cal U\big)\big]$ in $E_\alpha$ equals 
$\big(\Gamma^1_{{}_{f, L}}(S_\alpha)\cdot V(f)\big)_\bold 0 - 
\big(\Gamma^1_{{}_{f, L}}(S_\alpha)\cdot V(L)\big)_\bold 0$ if $S_\alpha\not\subseteq V(f)$. However, this is exactly the
content of Lemma 2.7.\qed

\vskip .4in

\noindent{\it Remark 2.11}.  In special cases, Theorem 2.10 was already known.

\vskip .1in

Consider the case where $X=\Cal U$ and
$\Fdot$ is the constant sheaf.  Then, $\operatorname{Ch}(\Fdot) =
\Cal U\times\{\bold 0\}$, and the image of $d\tilde f$ in $\Cal U\times\{\bold 0\}$ is simply defined by the Jacobian
ideal of $f$.  Hence, our result reduces to the result obtained from the work of Kashiwara in [{\bf K}] and
L\^e-Mebkhout in [{\bf L-M}] -- namely, that the projectivized characteristic cycle of the sheaf of vanishing cycles is
isomorphic to the exceptional divisor of the blow-up of the Jacobian ideal in affine space.

\vskip .1in

As a second special case, suppose that $X$ and $\Fdot$ are completely general, but that $\bold x$ is an isolated point
in the image of $\operatorname{Ch}(\phi_f\Fdot)$ in $X$ (for instance, $\bold x$ might be an isolated point in
$\operatorname{supp} \phi_f\Fdot$).  Then, for every stratum for which $m_\alpha\neq 0$,
$(\bold x, d_\bold x\tilde f)$ is an isolated point of $\imdf\cap \con$ or is not contained in the
intersection at all.  Now, $\con$ is an $(n+1)$-dimensional analytic variety and $\imdf$ is defined
by $n+1$ equations.  Therefore, $(\bold x, d_\bold x\tilde f)$ is regularly embedded in $\con$.

It follows that the exceptional divisor of the blow-up of $\imdf$ in $\con$ has one component over
$(\bold x, d_\bold x\tilde f)$ and that that component occurs with multiplicity precisely equal to the intersection
multiplicity
$\left(\imdf\cdot
\con\right)_{(\bold x, d_\bold x\tilde f)}$ in $T^*\Cal U$.  Thus, we recover the results of three independent works
appearing in [{\bf G}], [{\bf L\^e}], and [{\bf S}] -- that the coefficient of
$\{\bold x\}\times \Bbb C^{n+1}$ in $\operatorname{Ch}(\phi_f\Fdot)$ is given by 
$\left(\imdf\cdot \operatorname{Ch}(\Fdot)\right)_{(\bold x, d_\bold x\tilde f)}$.

\vskip .1in

In addition to generalizing the above results, Theorem 2.10 fits in well with Theorem 3.4.2 of [{\bf BMM}]; that theorem
contains a nice description of the characteristic cycles of the nearby cycles and of the restriction of a complex to a
hypersurface. However, [{\bf BMM}] does not contain a nice description of the vanishing cycles, nor does our Theorem
2.10 seem to follow easily from the results of [{\bf BMM}]; in fact, Example 3.4.3 of [{\bf BMM}] makes it clear that the general
result contained in our Theorem 2.10 was unknown -- for Brian\c con, Maisonobe, and Merle only derive the vanishing cycle result
from their nearby cycle result in the easy, known case where the
vanishing cycles are supported on an isolated point and, even then, they must make half a page of argument.

\vskip .4in

\noindent{\bf Corollary 2.12}. {\it For each extension $\tilde f$ of $f$, let $E_{\tilde f}$ denote the exceptional
divisor in
$\operatorname{Bl}_{\imdf}\overline{T^*_{{}_{X_{\operatorname{reg}}}}\Cal U}$. Then,
$\pi\big(E_{\tilde f}\big)$ is independent of $\tilde f$.}

\vskip .3in

\noindent{\it Proof}. We apply Theorem 2.10 to a complex of
sheaves $\Fdot$ such that
$m_\alpha =1$ for each smooth component of 
$X_{\operatorname{reg}}$ and $m_\alpha = 0$ for every other stratum in some Whitney stratification of $X$ (it is easy to
produce such an
$\Fdot$ -- see, for instance, Lemma 3.1 of [{\bf M1}]). The corollary follows from the fact that $\Bbb
P(\operatorname{Ch}(\phi_f\Fdot))$ does not depend on the extension. \qed

\vskip .5in

\noindent{\bigbold Section 3}. {\bf THE SPECIAL CASE OF PERVERSE SHEAVES}. 

\vskip .2in

We continue with our previous notation. 

\vskip .1in

For the purposes of this paper, perverse sheaves are important because the vanishing cycles functor  (shifted by
$-1$) applied to a perverse sheaf once again yields a perverse sheaf and because of the following lemma.

\vskip .2in

\noindent{\bf Lemma 3.1}. {\it  If $\Pdot$ is a perverse sheaf on $X$, then $\operatorname{Ch}(\Pdot)=\sum_\alpha
m_\alpha\conc$, where 
$$
m_\alpha = (-1)^{\dm X}\dm H^0(N_\alpha, \Bbb L_\alpha;\ \Pdot_{|_{N_\alpha}}[-d_\alpha]);
$$ 
in particular, 
$(-1)^{\dm X}\operatorname{Ch}(\Pdot)$ is a non-negative cycle.

If $\Pdot$ is perverse on $X$ (or, even, perverse up to a shift), then
$\operatorname{supp}\Pdot$ equals the image in $X$ of the characteristic cycle of $\Pdot$. }

\vskip .2in

\noindent{\it Proof}. The first statement follows from the definition of the characteristic cycle, together with the
fact that a perverse sheaf supported on a point has non-zero cohomology only in degree zero.

The second statement follows at once from the fact that if $\Pdot$ is perverse up to a shift, then so is the restriction
of $\Pdot$ to its support.  Hence, by the support condition on perverse sheaves, there is an open dense set of the
support, $\Omega$, such that, for all $\bold x\in \Omega$,
$H^*(\Pdot)_\bold x$ is non-zero in a single degree.  The conclusion follows. \qed

\vskip .2in

The fact that the above lemma refers to the support of $\Pdot$, which is the closure of the set of points with non-zero
stalk cohomology, means that we can use it to conclude something about the closure of the $\Pdot$-critical locus (recall
Definition 1.9).

\vskip .2in

\noindent{\bf Theorem 3.2}. {\it  Let $\Pdot$ be a perverse sheaf on $X$, and suppose that the characteristic cycle of
$\Pdot$ in $\Cal U$ is given by
$\operatorname{Ch}(\Pdot) = \sum_\alpha m_\alpha \conc$.

 Then, the closure of the
$\Pdot$-critical locus of
$f$ is given by
$$
\overline{\Sigma_{{}_{\Pdot}}f} \ =\ \Big\{\bold x\in X\ \big|\ (\bold x, d_\bold x\tilde
f)\in|\operatorname{Ch}(\Pdot)|\Big\}\ =\ \bigcup_{m_\alpha\neq 0}
\Sigma_{\operatorname{cnr}}\big(f_{|_{\overline{S}_\alpha}}\big).
$$ }

\vskip .2in

\noindent{\it Proof}. Let $\bold q\in X$, and let $v=f(\bold q)$. Let $\Cal W$ be an open neighborhood of $\bold q$ in
$X$ such that
$\Cal W\cap {\Sigma_{{}_{\Pdot}}f}\subseteq V(f-v)$ (see the end of Remark 1.7).   Then, $\Cal
W\cap\overline{\Sigma_{{}_{\Pdot}}f} =\Cal W\cap\operatorname{supp}\phi_{f-v}\Pdot$. As
$\phi_{f-v}\Pdot[-1]$ is perverse, Lemma 3.1 tells us that $\operatorname{supp}\phi_{f-v}\Pdot$ equals the image in
$X$ of
$\operatorname{Ch}(\phi_{f-v}\Pdot)$. Now, Theorem 2.10 tells us that this image is precisely
$$\bigcup_{m_\alpha\neq 0}
\Big\{\bold x\in \overline{S_\alpha}\ |\ (\bold x, d_\bold x\tilde f)\in\con\Big\},$$ since there can be no cancellation
as all the non-zero $m_\alpha$ have the same sign. 

Therefore, we have the desired equality of sets in an open neighborhood of every point; the theorem follows.\qed

\vskip .2in

 We will use the  {\it perverse cohomology} of the shifted constant sheaf, $\Bbb
C^\bullet_X[k]$, in order to deal with non-l.c.i.'s; this perverse cohomology is denoted by 
${}^p\negmedspace H^0(\Bbb C^\bullet_X[k])$ (see [{\bf BBD}] or [{\bf K-S}]). Like the intersection cohomology complex, 
this sheaf has the property that it is the shifted constant sheaf on
 the smooth part of any component of $X$ with dimension equal to $\dm X$.  
\vskip .2in 

We now list some properties of the perverse cohomology and of vanishing cycles that we will need later.  The reader is
referred to [{\bf BBD}] and [{\bf K-S}].

\bigskip

The perverse cohomology functor on $X$, 
${}^p\negmedspace H^0$, is a functor from the derived category of bounded, constructible complexes on $X$ to the Abelian category
of perverse sheaves on $X$.

The functor 
${}^p\negmedspace H^0$, applied to a perverse sheaf $\Pdot$ is canonically isomorphic to $\Pdot$. In addition, a bounded,
constructible complex of sheaves $\Fdot$ is perverse if and only ${}^p\negmedspace H^0(\Fdot[k])=0$ for all $k\neq 0$. In
particular, if $X$ is an l.c.i., then ${}^p\negmedspace H^0(\Bbb C^\bullet_X[\dm X])\cong \Bbb C^\bullet_X[\dm X]$ and
${}^p\negmedspace H^0(\Bbb C^\bullet_X[k]) = 0$ if $k\neq \dm X$.

The functor 
${}^p\negmedspace H^0$ commutes with vanishing cycles with a shift of $-1$, nearby cycles with a shift of $-1$, and 
Verdier dualizing.  That is, there are natural isomorphisms
$${}^p\negmedspace H^0 \circ \phi_f[-1] \cong
\phi_f[-1] \circ {}^p\negmedspace H^0,\hskip .2in {}^p\negmedspace H^0 \circ \psi_f[-1] \cong
\psi_f[-1] \circ {}^p\negmedspace H^0, \hskip .1in\text{and }\Cal D\circ {}^p\negmedspace H^0\cong {}^p\negmedspace
H^0\circ\Cal D.$$ 

Let $\bold F^\bullet$  be a bounded  complex of sheaves on $X$ which is constructible with respect to a connected
Whitney stratification $\{S_\alpha\}$ of $X$. Let
$S_{\operatorname{max}}$ be a maximal stratum contained in the support of $\Fdot$, and let  $m=\dm
S_{\operatorname{max}}$. Then, $\left({}^p\negmedspace H^0(\Fdot)\right)_{|_{S_{\operatorname{max}}}}$ is isomorphic (in
the derived category) to the complex which has $\left(\bold H^{-m}(\Fdot)\right)_{|_{S_{\operatorname{max}}}}$ in degree
$-m$ and zero in all other degrees. 

In particular, $\operatorname{supp}\Fdot = \bigcup_i \operatorname{supp}{}^p\negmedspace H^0(\Fdot[i])$, and if $\Fdot$
is supported on an isolated point, $\bold q$, then
$H^0({}^p\negmedspace H^0 (\bold F^\bullet))_{\bold q}
\cong H^0(\bold F^\bullet)_{\bold q}.$ 

\vskip .4in

Throughout the remainder of this paper, we let $\kp$ denote the perverse
sheaf ${}^p\negmedspace H^0(\Bbb C_X^\bullet[k+1])$; it will be useful later to have a nice characterization of the characteristic
cycle of $\kp$. 

\vskip .4in

\noindent{\bf Proposition 3.3}. {\it  
The complex 
$\kp$ is a perverse sheaf on $X$ which is constructible with respect to $\Cal S$ and the
characteristic cycle $\operatorname{Ch}(\kp)$ is equal to
$$(-1)^{\dm X}\sum_\alpha b_{k+1-d_\alpha}(N_\alpha, \Bbb L_\alpha)\conc,$$ where $b_j$ denotes the
$j$-th (relative) Betti number. }

\vskip .2in

\noindent{\it Proof}. The constructibility claim follows from the fact that the constant sheaf itself is clearly
constructible with respect to any Whitney stratification. The remainder follows trivially from the definition of the
characteristic cycle, combined with two properties of
${}^p\negmedspace H^0$; namely, 
${}^p\negmedspace H^0$ commutes with
$\phi_f[-1]$, and
${}^p\negmedspace H^0$ applied to a complex which is supported at a point simply gives ordinary cohomology in degree
zero and zeroes in all other degrees. See [{\bf K-S}, 10.3].
\qed

\vskip .3in

\noindent{\it Remark 3.4}. As $N_\alpha$ is contractible, it is possible to give a characterization of 
$b_{k+1-d_\alpha}(N_\alpha, \Bbb L_\alpha)$ without referring to $N_\alpha$; the statement gets a little complicated,
however, since we have to worry about what happens near degree zero and because the link of a maximal stratum is empty.
However, if we slightly modify the usual definitions of reduced cohomology and the corresponding reduced Betti numbers,
then the statement becomes quite easy.

What we want is for the ``reduced'' cohomology  $\widetilde H^k(A; \Bbb C)$ to be the relative cohomology vector space
$H^{k+1}(B, A;
\Bbb C)$, where $B$ is a contractible set containing $A$, and we want $\tilde b_*()$ to be the Betti numbers of this
``reduced'' cohomology. Therefore, letting $b_k()$ denote the usual $k$-th Betti number, we define $\tilde b_*()$ by
$$\tilde b_k(A) =\cases 
b_k(A), &\text{if $k\neq 0$ and $A\neq\emptyset$}\\
b_0(A)-1,&\text{if $k= 0$ and $A\neq\emptyset$}\\
0,&\text{if $k\neq -1$ and $A=\emptyset$}\\
1,&\text{if $k= -1$ and $A=\emptyset$.}
\endcases
$$
Thus, $\tilde b_k(A)$ is the $k$-th Betti number of the reduced cohomology, provided that $A$ is not the empty set.

\vskip .1in

 We let $\widetilde H^k(A; \Bbb C)$ denote the vector space ${\Bbb C}^{\tilde b_k(A)}$.

\vskip .1in

With this notation, the expression $b_{k+1-d_\alpha}(N_\alpha, \Bbb L_\alpha)$, which appears in 3.3, is equal to
$\tilde b_{k-d_\alpha}(\Bbb L_\alpha)$. The special definition of $\tilde b_k()$ for the empty set implies that if
$S_\alpha$ is maximal, then 
$$b_{k+1-d_\alpha}(N_\alpha, \Bbb L_\alpha) =\cases 0, &\text{if $k+1\neq d_\alpha$}\\
1,&\text{if $k+1=d_\alpha$}.
\endcases
$$

\vskip .2in

Hence, if $\dsize\operatorname{Ch}(\kp)=\sum_\alpha m_\alpha\big(\kp\big)\conc$,
then 3.3 implies that: $H^*(\Bbb L_\alpha; \Bbb C)\cong H^*(point; \Bbb C)$ if and only if $ m_\alpha\big(\kp\big)=0$ for all $k$.

\vskip .5in

By combining 3.2 with 3.3 and 3.4, we can now give a result about $\Sigma_{\Bbb C}f$. First, though, it will be useful
to adopt the following terminology.

\vskip .3in

\noindent{\bf Definition 3.5}. We say that the stratum $S_\alpha$ is {\it visible} (or, {\it $\Bbb C$-visible}) if
$H^*(\Bbb L_\alpha; \Bbb C)\not\cong H^*(point; \Bbb C)$ (or, equivalently, if $H^*(N_\alpha, \Bbb L_\alpha; \Bbb
C)\neq 0$). Otherwise, the stratum is {\it invisible}.

The final line of Remark 3.4 tells us that a stratum is visible if and only if there exists an integer $k$ such that
$\conc$  appears with a non-zero coefficient in $\operatorname{Ch}\big(\kp\big)$.

\vskip .1in

Note that if $S_\alpha$ has an empty complex link (i.e., the stratum is maximal), then $S_\alpha$ is {\bf visible}.

\vskip .3in

\noindent{\bf Theorem 3.6}. {\it   Then, 
$$
\overline{\Sigma_{\Bbb C}f}\ =\ \bigcup_{k=-1}^{\dm X-1} \overline{\Sigma_{{}_{\kp}}f}\ =\
\bigcup_{\text{visible }S_\alpha}\Big\{\bold x\in
\overline{S_\alpha}\ |\ (\bold x, d_\bold x\tilde f)\in\con\Big\}\ =\ \bigcup_{\text{visible }S_\alpha}
\Sigma_{\operatorname{cnr}}\big(f_{|_{\overline{S}_\alpha}}\big).
$$

In particular, since all maximal strata are visible, $\Sigma_{\operatorname{cnr}}f\subseteq\overline{\Sigma_{\Bbb C}f}$
(as stated in Proposition 1.6). Moreover, if $\bold x$ is an isolated point of $\Sigma_{\Bbb C}f$,then, for all Whitney stratifications,
$\{R_\beta\}$, of $X$, the only possibly visible stratum which can be contained in $f^{-1}f(\bold x)$ is $\{\bold x\}$.
}

\vskip .2in

\noindent{\it Proof}. Recall that, for any complex $\Fdot$, $\operatorname{supp}\Fdot = \bigcup_k
\operatorname{supp}{}^p\negmedspace H^0(\Fdot[k])$. In addition, we claim that $\kp=0$ unless
$-1\leqslant k\leqslant
\dm X-1$. By Lemma 3.1, $\kp=0$ is equivalent to $\operatorname{Ch}(\kp)=0$; if
$k$ is not between $-1$ and $\dm X-1$, then, using Proposition 3.3, $\operatorname{Ch}(\kp)=0$
follows from the fact that the complex link of a stratum has the homotopy-type of a finite CW complex of dimension no
more than the complex dimension of the link (see [{\bf Go-Mac1}]).

\vskip .1in

Now, in an open neighborhood of any point $\bold q$ with $v:= f(\bold q)$, we have
$$\overline{\Sigma_{\Bbb C}f} \ =\ \operatorname{supp}\phi_{f-v}\Bbb C^\bullet\ =\ \bigcup_k
\operatorname{supp}{}^p\negmedspace H^0(\phi_{f-v}\Bbb C_X^\bullet[k])\ = 
$$
$$
\bigcup_k \operatorname{supp}\phi_{f-v}[-1]\big({}^p\negmedspace H^0(\Bbb C_X^\bullet[k+1])\big)\ =\ \bigcup_k
\overline{\Sigma_{{}_{\kp}}f}.
$$
 Now, applying Theorem 3.2, we have
$$
\overline{\Sigma_{\Bbb C}f} \ =\ \bigcup_k\bigcup_{m_\alpha(\kp)\neq 0 }\Big\{\bold x\in
\overline{S_\alpha}\ |\ (\bold x, d_\bold x\tilde f)\in\con\Big\}.
$$ The desired conclusion follows.
\qed

\vskip .3in

\noindent{\it Remark 3.7}. Those familiar with stratified Morse theory should find the result of Theorem 3.6 very
un-surprising -- it looks like it results from some break-down of the $\Bbb C$-critical locus into normal and tangential
data, and naturally one gets no contributions from strata with trivial normal data. This is the approach that we took in
Theorem 3.2 of [{\bf Ma1}]. There is a slightly subtle, technical point which prevents us from taking this approach in
our current setting: by taking normal slices at points in an open, dense subset of $\operatorname{supp}\phi_{f-v}\Bbb
C^\bullet_X$, we could reduce ourselves to the case where $\Sigma_{\Bbb C}f$ consists of a single point, but we would
{\bf not}  know that the point was a {\bf stratified} isolated critical point. In particular, the case where
$\operatorname{supp}\phi_{f-v}\Bbb C^\bullet_X$ consists of a single point, but where $f$ has a non-isolated
(stratified) critical locus coming from an invisible stratum causes difficulties with the obvious Morse Theory approach.

\vskip .3in

\noindent{\it Remark 3.8}. At this point, we wish to add to our hierarchy of critical loci from Proposition 1.6.
Theorem 3.6 tells us that $\overline{\Sigma_{{}_{\kp}}f }\subseteq\overline{\Sigma_{{}_{\Bbb C}}f}$
for all $k$. If $X$ is purely $(d+1)$-dimensional, then 3.2 implies
that
$\Sigma_{\operatorname{cnr}}f\subseteq \overline{\Sigma_{{}_{{}^{d}\hskip -.01in\Pdot}}f}$.

Now, suppose that $X$ is irreducible of dimension $d+1$. Let $\text{\bf IC}^\bullet$  be the intersection
cohomology sheaf (with constant coefficients) on $X$ (see [{\bf Go-Mac2}]); $\text{\bf IC}^\bullet$ is a simple object in the
category of perverse sheaves. As the category of perverse sheaves on
$X$ is (locally) Artinian, and since 
${}^{d}\hskip -.01in\Pdot$ is a perverse sheaf which is the shifted constant sheaf on the smooth part of $X$, it
follows that $\text{\bf IC}^\bullet$ appears as a simple subquotient in any composition series for ${}^{d}\hskip
-.01in\Pdot$. Consequently, $|\operatorname{Ch}(\text{\bf IC}^\bullet)|\subseteq |\operatorname{Ch}({}^{d}\hskip
-.01in\Pdot)|$, and so 3.2 implies that $\overline{\Sigma_{{}_{\text{\bf IC}^\bullet}}f}\subseteq
\overline{\Sigma_{{}_{{}^{d}\hskip -.01in\Pdot}}f}$. Moreover, 3.2 also implies that
$\Sigma_{\operatorname{cnr}}f\subseteq \overline{\Sigma_{{}_{\text{\bf IC}^\bullet}}f}$. Therefore,  we can extend our
sequence of inclusions from Proposition 1.6 to:
$$
\overline{\Sigma_{\operatorname{reg}}f}\ \subseteq\
\overline{\Sigma_{\operatorname{alg}}f}\ \subseteq\ \overline{\Sigma_{{}_{\operatorname{Nash}}}f}\ \subseteq\
\Sigma_{\operatorname{cnr}}f\
\subseteq\
\overline{\Sigma_{{}_{\text{\bf IC}^\bullet}}f}\ \subseteq\ \overline{\Sigma_{{}_{{}^{d}\hskip -.01in\Pdot}}f}\subseteq\ 
\overline{\Sigma_{{}_\Bbb C}f}\ \subseteq\ \Sigma_{\operatorname{can}}f\ \subseteq\ \Sigma_{\operatorname{rdf}}f.
$$

Why not use one of these new critical loci as our most fundamental notion of the critical locus of $f$? Both
$\Sigma_{{}_{\text{\bf IC}^\bullet}}f$ and
$\Sigma_{{}_{{}^{d}\hskip -.01in\Pdot}}f$ are topological in nature, and easy examples show that they can be distinct
from $\Sigma_{{}_\Bbb C}f$. However, 3.6 tells us that
$\Sigma_{{}_{{}^{d}\hskip -.01in\Pdot}}f$ is merely one piece that goes into making up $\Sigma_{{}_{\Bbb C}}f$ -- we
should include the other shifted perverse cohomologies. On the other hand, given the importance of intersection
cohomology throughout mathematics, one should wonder why we do not use
$\Sigma_{{}_{\text{\bf IC}^\bullet}}f$ as our most basic notion. 

 Consider the node $X:= V(y^2-x^3-x^2)\subseteq\Bbb C^2$ and the function $f:=y_{|_X}$. The node has a {\it small
resolution of singularities} (see [{\bf Go-Mac2}]) given by simply pulling the branches apart. As a result, the
intersection cohomology sheaf on $X$ is the constant sheaf shifted by one on $X-\{\bold 0\}$, and the stalk
cohomology at $\bold 0$ is a copy of $\Bbb C^2$ concentrated in degree $-1$. Therefore, one can easily show that 
$\bold 0\not\in\overline{\Sigma_{{}_{\text{\bf IC}^\bullet}}f}$.

As $\overline{\Sigma_{{}_{\text{\bf IC}^\bullet}}f}$ fails to detect the simple change in topology of the level hypersurfaces
of $f$ as they go from being two points to being a single point, we do not wish to use $\Sigma_{{}_{\text{\bf IC}^\bullet}}f$
as our basic type of critical locus. That is not to say that $\Sigma_{{}_{\text{\bf IC}^\bullet}}f$ is not interesting in its
own right; it is integrally tied to resolutions of singularities. For instance, it is easy to show (using the
Decomposition Theorem [{\bf BBD}]) that if $\widetilde X@>\pi>> X$ is a resolution of singularities, then
$\Sigma_{{}_{\text{\bf IC}^\bullet}}f\subseteq
\pi(\Sigma(f\circ\pi))$. 

\vskip .5in

Now that we can ``calculate'' $\Sigma_{{}_{\Bbb C}}f$ using Theorem 3.6,  we
are ready to generalize the Milnor number of a function with an isolated critical point.

\vskip .3in

\noindent{\bf Definition 3.9}. If $\Pdot$ is a perverse sheaf on $X$, and $\bold x$ is an isolated point in
$\Sigma_{{}_{\Pdot}}f$ (or, if $\bold x\not\in\Sigma_{{}_{\Pdot}}f$), then we call
$\dm_{\Bbb C}H^0(\phi_{f-f(\bold x)}[-1]\Pdot)_\bold x$ the {\it Milnor number of $f$ at
$\bold x$ with coefficients in $\Pdot$} and we denote it by
$\mu_\bold x(f; \Pdot)$.

\vskip .1in

This definition is reasonable for, in this case, $\phi_{f-f(\bold x)}[-1]\Pdot$ is a perverse sheaf supported at the
isolated point
$\bold x$.  Hence, the stalk cohomology of 
$\phi_{f-f(\bold x)}[-1]\Pdot$ at $\bold x$ is possibly non-zero only in degree zero. Normally,  we summarize that
$\bold x$ is an isolated point in $\Sigma_{{}_{\Pdot}}f$ or that $\bold x\not\in\Sigma_{{}_{\Pdot}}f$  by writing
${\operatorname{dim}}_\bold x\Sigma_{{}_{\Pdot}}f\leqslant 0$ (we consider the dimension of the empty set to be
$-\infty$).

\vskip .4in

Before we state the next proposition, note that it is always the case that $\left(\imdf\ \cdot\ T^*_{\{\bold 0\}}\Cal U\right)_{(\bold 0,
d_\bold 0\tilde f)}=1$.

\vskip .4in

\noindent{\bf Proposition 3.10}. {\it For notational convenience, we assume that $\bold 0\in X$ and that
$f(\bold 0)=0$.

Then, ${\operatorname{dim}}_{\bold 0}\Sigma_{\Bbb C}f\leqslant 0$ if and only if,  for all $k$, 
${\operatorname{dim}}_{\bold 0}\Sigma_{\kp}f\leqslant 0$. Moreover, if ${\operatorname{dim}}_{\bold
0}\Sigma_{\Bbb C}f\leqslant 0$, then, 

\vskip .2in

\noindent i)\hskip .2in for all visible strata, $S_\alpha$, such that $\operatorname{dim}S_\alpha\geqslant
1$, the intersection of
$\imdf$ and
$\con$ is, at most,
$0$-dimensional at $(\bold 0, d_\bold 0\tilde f)$,
 
 and 
$$
\left(\imdf\ \cdot\ \con\right)_{(\bold 0, d_\bold 0\tilde f)}\ =\ \big(\Gamma^1_{{}_{f,
L}}(S_\alpha)\cdot V(f)\big)_\bold 0 - 
\big(\Gamma^1_{{}_{f, L}}(S_\alpha)\cdot V(L)\big)_\bold 0,
$$ where $L$ is a generic linear form, and

\vskip .2in

\noindent ii)\hskip .13in for all $k$, 
$$
\mu_\bold 0(f; \kp)\ =\ \tilde b_k(F_{f, \bold 0})\ =\ 
(-1)^{\operatorname{dim}X}\big(\imdf\ \cdot\ \operatorname{Ch}(\kp)\big)_{(\bold
0, d_\bold 0\tilde f)}\ =\ 
$$

$$\sum_{\text{visible }S_\alpha}\tilde b_{k-d_\alpha}(\Bbb L_\alpha)\left(\imdf\ \cdot\
\con\right)_{(\bold 0, d_\bold 0\tilde f)} =
$$

$$
\sum\Sb\text{visible }S_\alpha\\S_\alpha\text{ not maximal}\endSb \tilde b_{k-d_\alpha}(\Bbb
L_\alpha)\left(\imdf\ \cdot\ \con\right)_{(\bold 0, d_\bold 0\tilde f)} \ +\ \sum\Sb S_\alpha\text{
maximal}\\ \dm S_\alpha = k+1\endSb\left(\imdf\ \cdot\ \con\right)_{(\bold 0, d_\bold 0\tilde f)}. 
$$ 
 }

\vskip .3in

\noindent{\it Proof}.  It follows immediately from 3.6 that ${\operatorname{dim}}_{\bold 0}\Sigma_{\Bbb C}f\leqslant 0$
if and only if,  for all $k$, 
${\operatorname{dim}}_{\bold 0}\Sigma_{\kp}f\leqslant 0$. 

\vskip .1in

i) follows immediately from Lemma 2.7 (combined with Remark 2.11).

\vskip .1in

It remains for us to prove ii). As in the proof of 3.6, we have  
$${}^p\negmedspace H^0(\phi_{f}\Bbb C_X^\bullet[k])\ =\ \phi_{f}[-1]\big({}^p\negmedspace H^0(\Bbb
C_X^\bullet[k+1])\big)  =\
\phi_{f}[-1]\kp.$$ It follows that 
$$\mu_\bold 0(f; \kp)=\dm_{\Bbb C}H^0(\phi_{f}[-1]\kp)_\bold 0=\dm_{\Bbb
C}H^0\big({}^p\negmedspace H^0(\phi_{f}\Bbb C_X^\bullet[k])\big)_\bold 0 =\dm_{\Bbb C}H^0\big(\phi_{f}\Bbb
C_X^\bullet[k]\big)_\bold 0,$$ where the last equality is a result of the fact that $\bold 0$ is an isolated point in
the support of $\phi_{f}\Bbb C_X^\bullet[k]$ (for properties of ${}^p\negmedspace H^0$, see the beginning of Section 2).
Therefore, 
$$\mu_\bold 0(f; \kp)=\dm_{\Bbb C}H^0\big(\phi_{f}\Bbb C_X^\bullet[k]\big)_\bold 0 = \dm_{\Bbb
C}H^k\big(\phi_{f}\Bbb C_X^\bullet\big)_\bold 0 =
\operatorname{dim}\widetilde H^k(F_{f, \bold 0};\ \Bbb C).$$

That we also have the equality
$$
\mu_\bold 0(f; \kp)  = (-1)^{\operatorname{dim}X}\big(\imdf\ \cdot\
\operatorname{Ch}(\kp)\big)_{(\bold 0, d_\bold 0\tilde f)}
$$ is precisely the content of Theorem 2.10, interpreted as in the last paragraph of Remark 2.11.

The remaining equalities in ii) follow from the description of $\operatorname{Ch}(\kp)$ given in
Proposition 3.3 and Remark 3.4.\qed

\vskip .4in

\noindent{\it Remark 3.11}. The formulas from 3.10 provide a topological/algebraic method for ``calculating'' the Betti
numbers of the Milnor fibre for isolated critical points on arbitrary spaces. It should not be surprising that the data
that one needs is not just the algebraic data -- coming from the polar curves and intersection numbers -- but also
includes topological data about the underlying space: one has to know the Betti numbers of the complex links of strata. 

\vskip .3in

\noindent{\it Example 3.12}. The most trivial, non-trivial case where one can apply 3.10 is the case where $X$ is an
irreducible local, complete intersection with an isolated singularity (that is, $X$ is an irreducible i.c.i.s). Let us
assume that $\bold 0\in X$ is the only singular point of $X$ and that
$f$ has an isolated $\Bbb C$-critical point at $\bold 0$. Let $d$ denote the dimension of $X$.

Let us write
$\Bbb L_{X, \bold 0}$ for the complex link of $X$ at $\bold 0$. By [{\bf L\^e1}], $\Bbb L_{X, \bold 0}$  has the
homotopy-type of a finite bouquet of $(d-1)$-spheres.  Applying 3.10.ii, we see, then, that the reduced cohomology of
$F_{f,
\bold 0}$ is concentrated in degree $(d-1)$, and the $(d-1)$-th Betti number of  $F_{f, \bold 0}$ is equal to 
$$\tilde b_{d-1}(\Bbb L_{X, \bold 0})\left(\imdf\ \cdot\ T^*_\bold 0\Cal U\right)_{(\bold 0,
d_\bold 0\tilde f)} \ +\ \left(\imdf\ \cdot\
\overline{T^*_{X_{\operatorname{reg}}}\Cal U}\right)_{(\bold 0, d_\bold 0\tilde f)}=
$$
$$
\tilde b_{d-1}(\Bbb L_{X, \bold 0})\ +\ \big(\Gamma^1_{{}_{f, L}}(X_{\operatorname{reg}})\cdot V(f)\big)_\bold 0 - 
\big(\Gamma^1_{{}_{f, L}}(X_{\operatorname{reg}})\cdot V(L)\big)_\bold 0,
$$ for generic linear $L$.

Now, the polar curve and the intersection numbers are quite calculable in practice; see Remark 1.8 and Example 1.9 of
[{\bf Ma1}]. However, there remains the question of how one can compute $\tilde b_{d-1}(\Bbb L_{X, \bold 0})$. Corollary
4.6 and Example 5.4 of [{\bf Ma1}] provide an inductive method for computing the {\bf Euler characteristic} of $\Bbb
L_{X, \bold 0}$ (the induction is on the  codimension of $X$ in $\Cal U$) and, since we know that $\Bbb L_{X, \bold 0}$
has the homotopy-type of a bouquet of spheres, knowing the Euler characteristic is equivalent to knowing $\tilde
b_{d-1}(\Bbb L_{X, \bold 0})$.

The obstruction to using 3.10 to calculate Betti numbers in the general case is that, if $X$ is not an l.c.i., then a
formula for the Euler characteristic of the link of a stratum does not tell us the Betti numbers of the link.

\vskip .5in

\noindent{\bigbold Section 4}. {\bf THOM'S $a_f$ CONDITION}. 

\vskip .2in

We continue with the notation from Section 2.

\vskip .1in

In this section, we explain the fundamental relationship between Thom's $a_f$ condition and the vanishing cycles of $f$.

\vskip .2in

\noindent{\bf Definition 4.1}.  Let $M$ and $N$ be analytic submanifolds of $X$ such that
$f$ is constant on $N$. Then, the pair $(M, N)$ {\it satisfies Thom's $a_f$
condition at a point $\bold x\in N$} if and only if we have the containment
$\left(\overline{ T^*_{f_{|_M}}\Cal U}\right)_\bold x \ \subseteq \ 
\Big(T^*_{{}_N}\Cal U\Big)_\bold x$ of fibres over
$\bold x$.

\vskip .3in

We have been slightly more general in the above definition than is sometimes the case; we have not required that the rank of $f$ be constant on
$M$. Thus, if $X$ is an analytic space, we may write that $(X_{\operatorname{reg}}, N)$ satisfies the $a_f$ condition, instead of writing the much
more cumbersome $(X_{\operatorname{reg}} -
\Sigma\big(f_{|_{X_{\operatorname{reg}}}}\big), N)$ satisfies the $a_f$ condition. If $f$ is not constant on any irreducible component of $X$, it
is easy to see that these statements are equivalent:

Let $\overset\circ\to X := X_{\operatorname{reg}} -
\Sigma\big(f_{|_{X_{\operatorname{reg}}}}\big)$, which is dense in  $X_{\operatorname{reg}}$ (as $f$ is not constant on any irreducible components
of $X$).  We claim that
$\Bbb P\big(\overline{T^*_{f_{|_{\overset\circ\to X}}}\Cal U}\big) =
\Bbb P\big(\overline{T^*_{f_{|_{X_{\operatorname{reg}}}}}\Cal U}\big)$; clearly, this is equivalent to showing that
$T^*_{f_{|_{X_{\operatorname{reg}}}}}\Cal U \ \subseteq \ 
\overline{T^*_{f_{|_{\overset\circ\to X}}}\Cal U}$. This is simple, for if  $\bold x\in
\Sigma\big(f_{|_{X_{\operatorname{reg}}}}\big)$, then $(\bold x, \eta) \in  T^*_{f_{|_{X_{\operatorname{reg}}}}}\Cal U$
if and only if $(\bold x, \eta)\in T^*_{{}_{X_{\operatorname{reg}}}}\Cal U$, and $T^*_{{}_{X_{\operatorname{reg}}}}\Cal
U \ \subseteq \ 
\overline{T^*_{\overset\circ\to X}\Cal U} \ \subseteq \ \overline{T^*_{f_{|_{\overset\circ\to X}}}\Cal U}$. 

\vskip .5in

The link between Theorem 2.10 and the $a_f$ condition is provided by the following theorem, which describes the fibre in
the relative conormal in terms of the exceptional divisor in the blow-up of $\imdf$.  Originally,
we needed to assume Whitney's condition a) as an extra hypothesis; however, T. Gaffney showed us how to remove this
assumption by using a re-parameterization trick.

\vskip .2in

\noindent{\bf Theorem 4.2}. {\it  Let $\pi:\Cal U\times\Bbb C^{n+1}\times\Bbb P^n\rightarrow \Cal U\times\Bbb P^n$ 
denote the projection. 

Suppose that $f$ does not vanish identically on any irreducible component of $X$.   Let $E$ denote the exceptional
divisor in
$\operatorname{Bl}_{\imdf}\overline{T^*_{{}_{X_{\operatorname{reg}}}}\Cal U}\ \subseteq\ \Cal
U\times\Bbb C^{n+1}\times\Bbb P^n$.

Then, for all $\bold x\in X$, there is an inclusion of fibres over $\bold x$ given by $\big(\pi(E)\big)_\bold
x\subseteq\Big(\Bbb P\big(\overline{T^*_{f_{|_{X_{\operatorname{reg}}}}}\Cal U}\big)\Big)_\bold x$.  Moreover, if  $\bold
x\in\Sigma_{{}_{\operatorname{Nash}}}f$, then this inclusion is actually an equality.
 }

\vskip .3in

\noindent{\it Proof}.  By 2.12, it does not matter what extension of $f$ we use.  

\vskip .1in

That $\big(\pi(E)\big)_\bold x\subseteq \Big(\Bbb P\big(\overline{T^*_{f_{|_{X_{\operatorname{reg}}}}}\Cal
U}\big)\Big)_\bold x$ is easy. Suppose that $(\bold x, [\eta])\in \pi(E)$, that is $(\bold x, d_\bold x\tilde f,
[\eta])\in E$.  Then, there exists a sequence $(\bold x_i, \omega_i)\in \overline{T^*_{{}_{X_{\operatorname{reg}}}}\Cal
U}-\imdf$ such that $(\bold x_i, \omega_i, [\omega_i-d_{\bold x_i}\tilde f])\rightarrow (\bold x,
d_\bold x\tilde f, [\eta])$.  Hence, there exist scalars $a_i$ such that $a_i(\omega_i-d_{\bold x_i}\tilde f)\rightarrow
\eta$, and these $a_i(\omega_i-d_{\bold x_i}\tilde f)$ are relative conormal covectors whose projective class approaches
that of $\eta$. Thus, $\big(\pi(E)\big)_\bold x\subseteq \Big(\Bbb
P\big(\overline{T^*_{f_{|_{X_{\operatorname{reg}}}}}\Cal U}\big)\Big)_\bold x$. 

\vskip .2in

 We must now show that $\Big(\Bbb P\big(\overline{T^*_{f_{|_{X_{\operatorname{reg}}}}}\Cal U}\big)\Big)_\bold x
 \subseteq\big(\pi(E)\big)_\bold x$, provided that $\bold
x\in\Sigma_{{}_{\operatorname{Nash}}}f$.  

\vskip .1in  

Let $\overset\circ\to X := X_{\operatorname{reg}} - \Sigma\big(f_{|_{X_{\operatorname{reg}}}}\big)$.
 Suppose  that $(\bold x, [\eta])\in 
\Bbb P\big(\overline{T^*_{f_{|_{\overset\circ\to X}}}\Cal U}\big)$. Then, there exists a complex analytic
path $\alpha(t) = (\bold x(t), \eta_t)\in \overline{T^*_{f_{|_{\overset\circ\to X}}}\Cal U}$ such that
$\alpha(0) = (\bold x, \eta)$ and $\alpha(t)\in T^*_{f_{|_{\overset\circ\to X}}}\Cal U$ for $t\neq 0$. As
$f$ has no critical points on $\overset\circ\to X$, each $\eta_t$ can be written uniquely as $\eta_t =
\omega_t + \lambda(\bold x(t))d_{\bold x(t)}\tilde f$, where $\omega_t(T_{\bold x(t)}\overset\circ\to X) = 0$ and
$\lambda(\bold x(t))$ is a scalar. By evaluating each side on $\bold x^\prime(t)$, we find that  $\lambda(\bold x(t)) =
\frac{\eta_t(\bold x^\prime(t))}{\ \frac{d}{dt}f(\bold x(t))\ }$. 

Thus, as $\lambda(\bold x(t))$ is a quotient of two
analytic functions, there are only two possibilities for what happens to
$\lambda(\bold x(t))$ as $t\rightarrow 0$.

\vskip .3in

\noindent {\bf Case 1}:  $|\lambda(\bold x(t))|\rightarrow\infty$ as $t\rightarrow 0$.

\vskip .2in

In this case, since $\eta_t\rightarrow\eta$, it follows that 
$\dsize\frac{\eta_t}{\lambda(\bold x(t))}\rightarrow 0$ and, hence, $\dsize-\frac{\omega_t}{\lambda(\bold
x(t))}\rightarrow d_\bold x\tilde f$.   Therefore, $$\left(\bold x(t), -\frac{\omega_t}{\lambda(\bold x(t))},
\left[-\frac{\omega_t}{\lambda(\bold x(t))}-d_{\bold x(t)}\tilde f
\right]\right) \ = \ \left(\bold x(t), -\frac{\omega_t}{\lambda(\bold x(t))}, \left[\eta_t(\bold
x(t))\right]\right)\rightarrow (\bold x,
d_\bold x\tilde f, [\eta]),$$ and so $(\bold x, [\eta])\in \pi(E)$.
 
\vskip .4in

\noindent {\bf Case 2}:  $\lambda(\bold x(t))\rightarrow\lambda_0$ as $t\rightarrow 0$.

\vskip .2in

In this case, $\omega_t$ must possess a limit as $t\rightarrow 0$. For $t$ small and unequal to zero, let $\operatorname{proj}_t$
denote the complex orthogonal projection from the fibre 
$\big (T^*_{f_{|_{\overset\circ\to X}}}\Cal U\big)_{\bold x(t)}$ to the fibre $\big (T^*_{\overset\circ\to X}\Cal U\big)_{\bold
x(t)}$. Let $\gamma_t:=\operatorname{proj}_t(\eta_t)=\omega_t+\lambda(\bold x(t))\operatorname{proj}_t(d_{\bold x(t)}\tilde f)$. 
Since
$\bold x\in\Sigma_{{}_{\operatorname{Nash}}}f$, we have that $\operatorname{proj}_t(d_{\bold x(t)}\tilde f)\rightarrow d_\bold x\tilde
f$ and, thus, $\gamma_t\rightarrow\eta$.

As $\eta$ is not zero (since it represents a projective class), we may define the (real, non-negative) scalar 
$$
a_t:=\sqrt{\frac{||\operatorname{proj}_t(d_{\bold x(t)}\tilde f)-d_{\bold x(t)}\tilde f||}{||\gamma_t||}}.
$$
One now verifies easily that
$$
(\bold x(t), \ a_t\gamma_t+ \operatorname{proj}_t(d_{\bold x(t)}\tilde f), \ [a_t\gamma_t+ \operatorname{proj}_t(d_{\bold x(t)}\tilde
f)-d_{\bold x(t)}\tilde f])\ \longrightarrow\ (\bold x, d_\bold x\tilde f, [\eta]),
$$
and, hence, that $(\bold x, [\eta])\in \pi(E)$.
 \qed

\vskip .5in

\noindent{\it Remark 4.3}. In a number of results throughout the remainder of this paper, the reader will find the hypotheses that
 $\bold x\in \Sigma_{{}_{\operatorname{Nash}}}f$ or that $\bold x\in \Sigma_{\operatorname{alg}}f$.  While Theorem 4.2 explains why
the hypothesis $\bold x\in \Sigma_{{}_{\operatorname{Nash}}}f$ is important, it is not so clear why the hypothesis $\bold x\in
\Sigma_{\operatorname{alg}}f$ is of interest.

If $Y$ is an analytic subset of $X$, then one shows easily that $Y\cap\Sigma_{\operatorname{alg}}f\subseteq 
\Sigma_{\operatorname{alg}}(f_{|_Y})$. The Nash critical does not possess such an inheritance property. Thus, the easiest
hypothesis to make in order to guarantee that a point, $\bold x$, is in the Nash critical locus of any analytic subset containing
$\bold x$ is the hypothesis that $\bold x\in\Sigma_{\operatorname{alg}}f$, for then if $\bold x\in Y$, we conclude that $\bold x\in
\Sigma_{\operatorname{alg}}(f_{|_Y})\subseteq \Sigma_{{}_{\operatorname{Nash}}}(f_{|_Y})$.

\vskip .5in

We come now to the result which tells one how the topological information provided by the sheaf of vanishing cycles
controls the $a_f$ condition.

\vskip .3in

\noindent{\bf Corollary 4.4}. {\it Let $N$ be a submanifold of $X$ such that $N\subseteq V(f)$, and let $\bold x\in N$

Let $\operatorname{Ch}(\Fdot) = \sum_\alpha m_\alpha\big[\overline{T^*_{{}_{M_\alpha}}\Cal U}\big]$, where
$\{M_\alpha\}$ is a collection of connected analytic submanifolds of
$X$ such that either
$m_\alpha\geqslant 0$ for all $\alpha$, or $m_\alpha\leqslant 0$ for all $\alpha$.   Let
$\operatorname{Ch}(\phi_f\Fdot) = \sum_\beta k_\beta
\big[\overline{T^*_{{}_{R_\beta}}\Cal U}\big]$, where $\{R_\beta\}$ is a collection of connected analytic  submanifolds. 

Finally, suppose that, for all $\beta$, there is an inclusion of fibres over $\bold x$ given by  
$\Big(\overline{T^*_{{}_{R_\beta}}\Cal U}\Big)_\bold x \ \subseteq \ 
\left(T^*_{{}_N}\Cal U\right)_\bold x$.

\vskip .1in

Then, the pair
$\left(\left(\overline{M_\alpha}\right)_{\operatorname{reg}},
N\right)$ satisfies Thom's
$a_f$ condition at $\bold x$ for every 
$M_\alpha$ for which $f_{|_{M_\alpha}}\not\equiv 0$, $m_\alpha\neq 0$ and such that $\bold
x\in\Sigma_{{}_{\operatorname{Nash}}}(f_{|_{{\overline{M}}_\alpha}})$. }

\vskip .3in

\noindent{\it Proof}.  Let $\big\{S_\gamma\big\}$ be a Whitney stratification for $X$ such that each
$\overline{M_\alpha}$ is a union of strata and such that $\Sigma\big(f_{|_{S_\gamma}}\big) = \emptyset$ unless
$S_\gamma\subseteq V(f)$.  Hence, for each $\alpha$, there exists a unique
$S_\gamma$ such that $\overline{M_\alpha} = \overline{S_\gamma}$; denote this stratum by $S_\alpha$.  It follows at once
that $\operatorname{Ch}(\Fdot) = \sum_\alpha m_\alpha\big[\overline{T^*_{{}_{S_\alpha}}\Cal U}\big]$.

\vskip .1in

From Theorem 2.10, $E = \sum_\alpha m_\alpha E_\alpha \cong \Bbb P\big(\operatorname{Ch}(\phi_f\Fdot)\big)$.  Thus,
since all non-zero $m_\alpha$ have the same sign, if $m_\alpha$ is not zero, then $E_\alpha$ appears with a non-zero
coefficient in $\Bbb P\big(\operatorname{Ch}(\phi_f\Fdot)\big)$.  

The result now follows immediately by applying Theorem 4.2 to each
$\overline{M_\alpha}$ in place of $X$.
\qed

\vskip .3in

Theorem 4.2 also allows us to prove an interesting relationship between the characteristic varieties of the vanishing
and nearby cycles -- provided that the complex of sheaves under consideration is perverse.

\vskip .3in

\noindent{\bf Corollary 4.5}. {\it Let $\Pdot$ be a perverse sheaf on $X$. If $\bold x\in\Sigma_{\operatorname{alg}}f$ 
and $(\bold x, \eta)\in|\operatorname{Ch}(\psi_f\Pdot)|$, then
$(\bold x,
\eta)\in|\operatorname{Ch}(\phi_f\Pdot)|$. }

\vskip .3in

\noindent{\it Proof}. Let $\Cal S:=\{S_\alpha\}$ be a Whitney stratification with connected strata such that $\Pdot$ is
constructible with respect to $\Cal S$ and such that $V(f)$ is a union of strata. For the remainder of the proof, we
will work in a neighborhood of $V(f)$ -- a neighborhood in which, if $S_\alpha\not\subseteq V(f)$, then
$\Sigma(f_{|_{S_\alpha}})=\emptyset$.

Let $\operatorname{Ch}(\Pdot) = \sum m_\alpha\conc$. As $\Pdot$ is perverse, all non-zero $m_\alpha$ have the same sign.
Thus, 2.10 tells us -- using the notation from 2.10 -- that $$\dsize|\Bbb
P(\operatorname{Ch}(\phi_f\Pdot))|=\bigcup_{m_\alpha\neq 0}\pi(E_\alpha),\tag{$\dagger$}$$ where
$E_\alpha$ denotes the exceptional divisor in the blow-up of $\con$ along $\imdf$ (in a
neighborhood of
$V(f)$). In addition, 2.3 tells us that 
$$|\operatorname{Ch}(\psi_f\Pdot)|=\big(V(f)\times\Bbb C^{n+1}\big)\ \cap\bigcup_{\Sb m_\alpha\neq
0\\S_\alpha\not\subseteq V(f)\endSb}\overline{T^*_{f_{|_{S_\alpha}}}\Cal U}.$$

Assume  $(\bold x, \eta)\in|\operatorname{Ch}(\psi_f\Pdot)|$.  Then, there exists
$S_\alpha\not\subseteq V(f)$ such that
$m_\alpha\neq 0$ and
$(\bold x, \eta)\in\overline{T^*_{f_{|_{S_\alpha}}}\Cal U}$. Clearly, then, $(\bold x,
\eta)\in\overline{T^*_{f_{|_{{(\overline{S}_\alpha)}_{\operatorname{reg}}}}}\Cal U}$. Now, if $\bold
x\in\Sigma_{\operatorname{alg}}f$ and $\eta\neq 0$, then $\bold
x\in\Sigma_{\operatorname{alg}}(f_{|_{\overline{S}_\alpha}})$ and so Theorem 4.2 implies that
$(\bold x, [\eta])\in\pi(E_\alpha)$, where $[\eta]$ denotes the projective class of $\eta$ and $E_\alpha$ denotes the
exceptional divisor of the blow-up of $\overline{T^*_{{(\overline{S}_\alpha)}_{\operatorname{reg}}}\Cal U} = \con$ along
$\imdf$. Thus, by ($\dagger$), $(\bold x, \eta)\in|\operatorname{Ch}(\phi_f\Pdot)|$.

We are left with the trivial case of when $(\bold x, 0)\in|\operatorname{Ch}(\psi_f\Pdot)|$. Note that, if $(\bold x,
0)\in|\operatorname{Ch}(\psi_f\Pdot)|$, then there must exist some non-zero $\eta$ such that $(\bold x,
\eta)\in|\operatorname{Ch}(\psi_f\Pdot)|$. For, otherwise, the stratum (in some Whitney stratification) of
$\operatorname{supp}\psi_f\Pdot$ containing
$\bold x$ must be all of
$\Cal U$. However,
$\psi_f\Pdot$ is supported on $V(f)$, and so $f$ would have to be zero on all of $\Cal U$; but, this implies that
$|\operatorname{Ch}(\psi_f\Pdot)|=\emptyset$. Now, if we have some non-zero $\eta$ such that $(\bold x,
\eta)\in|\operatorname{Ch}(\psi_f\Pdot)|$, then by the above argument,  $(\bold x,
\eta)\in|\operatorname{Ch}(\phi_f\Pdot)|$ and, thus, certainly $(\bold x, 0)\in|\operatorname{Ch}(\phi_f\Pdot)|$.\qed 

\vskip .5in

\noindent{\bigbold Section 5}. {\bf CONTINUOUS FAMILIES OF CONSTRUCTIBLE COMPLEXES}. 

\vskip .2in

We wish to prove statements of the form: the constancy of certain data in a family implies that some nice geometric
facts hold. As the reader should have gathered from the last section, it is very advantageous to use complexes of
sheaves for cohomology coefficients; in particular, being able to use perverse coefficients is very desirable. The
question arises: what should a family of complexes mean?

\vskip .1in

Let $X$ be an analytic space, let $t:X\rightarrow \Bbb C$ be an analytic function, and let $\Fdot$ be a bounded,
constructible complex of
$\Bbb C$-vector spaces. We could say that $\Fdot$ and $t$ form a ``nice'' family of complexes, since, for all
$a\in\Bbb C$, we can consider the complex $\Fdot_{|_{t^{-1}(a)}}$ on the space $X_{|_{t^{-1}(a)}}$. However, this does
yield a satisfactory theory, because there may be absolutely no relation between $\Fdot_{|_{t^{-1}(0)}}$ and
$\Fdot_{|_{t^{-1}(a)}}$ for $a$ close to $0$. What we need is a notion of {\it continuous} families of complexes -- we
want $\Fdot_{|_{t^{-1}(0)}}$ to equal the ``limit'' of $\Fdot_{|_{t^{-1}(a)}}$ as $a$ approaches $0$. Fortunately, such
a notion already exists; it just is not normally thought of as continuity.

\vskip .3in

\noindent{\bf Definition 5.1}. Let $X$, $t$, and $\Fdot$ be as above. We define the {\it limit of
$\Fdot_a:=\Fdot_{|_{t^{-1}(a)}}[-1]$ as $a$ approaches $b$}, $\dsize\lim_{a\rightarrow b}\Fdot_a$, to be the nearby
cycles
$\psi_{t-b}\Fdot[-1]$.

We say that the family $\Fdot_a$ is {\it continuous at the value $b$} if the comparison map from $\Fdot_b$ to
$\psi_{t-b}\Fdot[-1]$ is an isomorphism, i.e., if the vanishing cycles $\phi_{t-b}\Fdot[-1] = 0$. We say that the family
$\Fdot_a$ is {\it continuous} if it is continuous for all values $b$.

We say that the family $\Fdot_a$ is {\it continuous at the point $\bold x\in X$} if there is an open neighborhood $\Cal
W$ of
$\bold x$ such that the family defined by restricting $\Fdot$ to $\Cal W$ is continuous at the value $t(\bold x)$. 

If $\Pdot$ is a perverse sheaf on $X$ and $\Pdot_a:=\Pdot_{|_{t^{-1}(a)}}[-1]$ is a continuous family of complexes, then
we say that $\Pdot_a$ is a continuous family of perverse sheaves.

\vskip .3in

\noindent{\bf Remark 5.2}. The reason for the shifts by $-1$ in the families is so that if $\Pdot$ is perverse,  and
$\Pdot_a$ is a continuous family, then each $\Pdot_a$ is, in fact, a  perverse sheaf (since
$\Pdot_a\cong\psi_{t-a}\Pdot[-1]$).

It is not difficult to show that: if the family $\Fdot_a$ is continuous at the value $b$, and, for all $a\neq b$, each
$\Fdot_a$ is perverse, then, near the vlaue $b$, the family $\Fdot_a$ is a continuous family of perverse sheaves.

\vskip .2in

For the remainder of this section, we will be using the following additional notation. Let $\tilde t$ be an analytic
function on
$\Cal U$, and let $t$ denote its restriction to $X$. Let
$\Pdot$ be a perverse sheaf on
$X$. Consider the families of spaces, functions, and sheaves given by $X_a:=X\cap V(t-a)$, $f_a:=f_{|_{X_a}}$, and
$\Pdot_a:=\Pdot_{|_{X_a}}[-1]$ (normally, if we are not looking at a specific value for $t$, we write $X_t$,
$f_t$, and
$\Pdot_t$ for these families). Note that, if we have as an hypothesis that
$\Pdot_t$ is continuous, then the family $\Pdot_t$ is actually a family of {\bf perverse} sheaves.

\vskip .4in

We will now prove three fundamental lemmas; all of them have trivial proofs, but they are nonetheless extremely useful. 

\vskip .1in

The first lemma uses Theorem 3.2 to characterize continuity at a point for families of perverse sheaves.

\vskip .3in

\noindent{\bf Lemma 5.3}. {\it Let $\bold x\in X$. The following are equivalent:

\vskip .1in

\noindent i)\hskip .2in The family $\Pdot_t$ is continuous at $\bold x$;

\vskip .2in

\noindent ii)\hskip .17in $\bold x\not\in\overline{\Sigma_{{}_{\Pdot}}t}$;

\vskip .2in

\noindent iii)\hskip .12in $(\bold x, d_{\bold x}\tilde t)\not\in |\operatorname{Ch}(\Pdot)|$ for some local extension,
$\tilde t$, of $t$ to
$\Cal U$ in a neighborhood of $\bold x$; and

\vskip .2in

\noindent iv)\hskip .16in $(\bold x, d_{\bold x}\tilde t)\not\in |\operatorname{Ch}(\Pdot)|$ for every local extension,
$\tilde t$, of $t$ to
$\Cal U$ in a neighborhood of $\bold x$. }

\vskip .3in

\noindent{\it Proof}. The equivalence of i) and  ii) follows from their definitions, together with Remark 1.7. The
equivalence between ii), iii), and iv) follows immediately from Theorem 3.2. 
\qed

\vskip .3in

The next lemma is a necessary step in several proofs.

\vskip .3in

\noindent{\bf Lemma 5.4}. {\it  Suppose that the family $\Pdot_t$ is continuous at $t=b$, and that the characteristic
cycle of
$\Pdot$ is given by $\sum_\alpha m_\alpha\conc$. Then, $S_\alpha\not\subseteq V(t-b)$ if $m_\alpha\neq 0$. }

\vskip .3in

\noindent{\it Proof}. This follows immediately from 5.3.
\qed

\vskip .4in

The last of our three lemmas is the  {\bf stability of continuity} result.

\vskip .3in

\noindent{\bf Lemma 5.5}. {\it  Suppose that the family $\Pdot_t$ is continuous at $\bold x\in X$. Then, $\Pdot_t$ is
continuous at all points near $\bold x$. In addition, if $\overset\circ\to{\Bbb D}$ is an open disk around the origin in
$\Bbb C$, $h: \overset\circ\to{\Bbb D}\times X\rightarrow
\Bbb C$ is an analytic function, $h_c(\bold z):=h(c, \bold z)$, and $h_0= t$, then the family $\Pdot_{h_c}$ is
continuous at $\bold x$ for all $c$ sufficiently close to $0$. }

\vskip .3in

\noindent{\it Proof}. Let $\tilde t$ be an extension of $t$ to a neighborhood of $\bold x$ in $\Cal U$, and let
$\Pi_1:T^*\Cal U\rightarrow\Cal U$ be the cotangent bundle. As
$T^*\Cal U$ is isomorphic to
$\Cal U\times
\Bbb C^{n+1}$, there is a second projection $\Pi_2:T^*\Cal U\rightarrow\Bbb C^{n+1}$. 

Now,
$\Pi_1^{-1}({\bold x})
\cap |\operatorname{Ch}(\Pdot)|$ and $\Pi_2^{-1}(d_\bold x\tilde t)
\cap |\operatorname{Ch}(\Pdot)|$ are closed sets. Therefore, the lemma follows immediately from 5.3.
\qed

\vskip .4in

The following lemma allows us to use intersection-theoretic arguments for families of generalized isolated critical
points. 

\vskip .3in

\noindent{\bf Lemma 5.6}. {\it   Suppose that the family $\Pdot_t$ is continuous at
$\bold x\in X$. Let $b:= t(\bold x)$. Let
$\{S_\alpha\}$ be a Whitney stratification of $X$ with connected strata with respect to which $\Pdot$ is constructible
and such that
$V(t-b)$ is a union of strata. Suppose that $\operatorname{Ch}(\Pdot)$ is given by $\sum_\alpha m_\alpha\conc$. If
$\dm_\bold x\Sigma_{{}_{\Pdot_b}}f_b\leqslant 0$, then there exists an open neighborhood $\Cal W$ of $\bold x$ in $\Cal
U$ such that:

\vskip .2in

\noindent i)\hskip .2in $\imdf$ properly intersects 
$\dsize\sum_{\alpha}m_\alpha\Big[\overline{T^*_{t_{|_{S_\alpha}}}\Cal U}\Big]$ in $\Cal W$;

\vskip .3in

\noindent ii)\hskip .16in for all $\bold y\in X\cap \Cal W$, $V(t-t(\bold y))$ properly intersects 
$$
\imdf\ \cdot\ \dsize\sum_{\alpha}m_\alpha\Big[\overline{T^*_{t_{|_{S_\alpha}}}\Cal U}\Big]
$$ at $(\bold y, d_\bold y\tilde f)$ in (at most) an isolated point; and 

\vskip .3in

\noindent iii)\hskip .12in for all $\bold y\in X\cap \Cal W$, if $a:= t(\bold y)$, then $\dm_\bold
y\Sigma_{{}_{\Pdot_a}}f_a\leqslant 0$ and 
$$
\mu_\bold y(f_a; \Pdot_a) = (-1)^{\dm X}\Big[\Big(\imdf\ \cdot\
\sum_{\alpha}m_\alpha\Big[\overline{T^*_{t_{|_{S_\alpha}}}\Cal U}\Big]\Big)\ \cdot\ V(t-a)\Big ]_{(\bold y, d_\bold
y\tilde f)}.
$$
  }

\vskip .3in

\noindent{\it Proof}.  First, note that we may assume that $X = \operatorname{supp}\Pdot$; for, otherwise, we would
immediately replace $X$ by
$\operatorname{supp}\Pdot$. Now, it follows from Lemma 5.4 that $V(t-b)$ does not contain an entire irreducible
component of $X$. Thus $\dm X_0 = \dm X-1$.

We use $\tilde f$ as a common extension of $f_t$ to $\Cal U$, for all $t$. Proposition 3.10 tells us that $\mu_\bold
x(f_b; \Pdot_b) = (-1)^{\dm X-1}\big(\imdf\ \cdot\ \operatorname{Ch}(\Pdot_b)\big)_{(\bold b,
d_\bold b\tilde f)}$.  Then, continuity, implies that $\operatorname{Ch}(\Pdot_b)
=\operatorname{Ch}(\psi_{t-b}[-1]\Pdot)$, and 
$$
\operatorname{Ch}(\psi_{t-b}[-1]\Pdot) =-\operatorname{Ch}(\psi_{t-b}\Pdot) = -\big(V(t-b)\times\Bbb C^{n+1}\big)\ \cdot\
\sum_{S_\alpha\not\subseteq V(t-b)}m_\alpha\Big[\overline{T^*_{t_{|_{S_\alpha}}}\Cal U}\Big],\tag{$*$}
$$  by Theorem 2.3. By Lemma 5.4, we may index over all $S_\alpha$; for, if $S_\alpha\subseteq V(t-b)$, then $m_\alpha =
0$.

Therefore,
$$
\mu_\bold x(f_b; \Pdot_b) = (-1)^{\dm X}\Big(\imdf\ \cdot\ \big(V(t-b)\times\Bbb C^{n+1}\big)\
\cdot\
\sum_{\alpha}m_\alpha\Big[\overline{T^*_{t_{|_{S_\alpha}}}\Cal U}\Big]\Big)_{(\bold x, d_\bold x\tilde
f)}=\tag{$\dagger$}
$$
$$ (-1)^{\dm X}\Big(\Big(\imdf\ \cdot\
\sum_{\alpha}m_\alpha\Big[\overline{T^*_{t_{|_{S_\alpha}}}\Cal U}\Big]\Big)\ \cdot\ \big(V(t-b)\times\Bbb
C^{n+1}\big)\Big)_{(\bold x, d_\bold x\tilde f)}.
$$ Thus, $$C:=(-1)^{\dm X}\Big(\imdf\ \cdot\
\sum_{\alpha}m_\alpha\Big[\overline{T^*_{t_{|_{S_\alpha}}}\Cal U}\Big]\Big)$$ is a non-negative cycle such that $(\bold
x, d_\bold x\tilde f)$ is an isolated point in (or, is not in) $C\cdot V(t-b)$. Statements i) and ii) of the lemma
follow immediately.

Now, Lemma 5.5 tells us that the family $\Pdot_t$ is continuous at all points near $\bold x$; therefore, if $\bold y$ is
close to
$\bold x$ and $a:=t(\bold y)$, then, by repeating the argument for $(*)$, we find that
$$
\operatorname{Ch}(\Pdot_a) =-\operatorname{Ch}(\psi_{t-a}\Pdot) = -\big(V(t-a)\times\Bbb C^{n+1}\big)\ \cdot\
\sum_{\alpha}m_\alpha\Big[\overline{T^*_{t_{|_{S_\alpha}}}\Cal U}\Big]
$$ and we know that the intersection of this cycle with
$\imdf$ is (at most) zero-dimensional at
$(\bold y, d_\bold y\tilde f)$ (since $C\cap V(t-b)$ is (at most) zero-dimensional at $\bold x$). By considering $\tilde
f$ an extension of $f_a$ and applying Theorem 3.2, we conclude that
$\dm_\bold y\Sigma_{{}_{\Pdot_a}}f_a\leqslant 0$. 

Finally, now that we know that $\Pdot_t$ is continuous at $\bold y$ and that $\dm_\bold
y\Sigma_{{}_{\Pdot_a}}f_a\leqslant 0$, we may argue as we did at $\bold x$ to conclude that $(\dagger)$ holds with
$\bold x$ replaced by $\bold y$ and $b$ replaced by $a$. This proves iii).
\qed

\vskip .4in

We can now prove an {\bf additivity/upper-semicontinuity} result. We prove this result for a more general type of family
of perverse sheaves; instead of parametrizing by the values of a function, we parametrize implicitly. We will need this
more general perspective in Theorem 5.10.

\vskip .3in

\noindent{\bf Theorem 5.7}. {\it Suppose that the family $\Pdot_t$ is continuous at
$\bold x\in X$. Let $b:= t(\bold x)$, and suppose that $\dm_\bold x\Sigma_{{}_{\Pdot_b}}f_b\leqslant 0$. 

Let $\overset\circ\to{\Bbb D}$ be an open disk around the origin in $\Bbb C$, let $h: \overset\circ\to{\Bbb D}\times
X\rightarrow
\Bbb C$ be an analytic function, for all $c\in\overset\circ\to{\Bbb D}$, let $h_c(\bold z):=h(c, \bold z)$,  let
${}_c\Pdot:= \Pdot_{|_{V(h_c-b)}}[-1]$ and ${}_cf := f_{|_{V(h_c-b)}}$. Suppose that $h_0= t$.

Then, there exists an open neighborhood $\Cal W$ of $\bold x$ in $\Cal U$ such that, for all small $c$, for all $\bold
y\in V(h_c-b)\cap\Cal W$, $\dm_\bold y\Sigma_{{}_{{}_c\Pdot}}{}_cf\leqslant 0$.

 Moreover, for fixed $c$ close to $0$, there are a finite number of points $\bold y\in V(h_c-b)\cap\Cal W$ such that 
$\mu_\bold y({}_cf; {}_c\Pdot)\neq 0$ and
$$
\mu_\bold x({}_bf; {}_b\Pdot) = \sum_{\bold y\in  V(h_c-b)\cap\Cal W}\mu_\bold y({}_cf; {}_c\Pdot). 
$$

\vskip .2in

In particular, for all small $c$, for all $\bold y\in V(h_c-b)\cap\Cal W$, $\mu_\bold y({}_cf; {}_c\Pdot)
\leqslant \mu_\bold x({}_bf; {}_b\Pdot)$. }

\vskip .3in

\noindent{\it Proof}. We continue to let $\Pdot_c=\Pdot_{|_{V(t-c)}}[-1]$ and $f_c=f_{|_{V(t-c)}}$. Note that, if we let
$h(w, \bold z):= t(\bold z)-w$, then the statement of the theorem would reduce to a statement about the ordinary
families  $\Pdot_c$ and $f_c$. Moreover, this statement about the families $\Pdot_c$ and $f_c$ follows immediately from
Lemma 5.6. We wish to see that this apparently weak form of the theorem actually implies the stronger form.

\vskip .1in

Shrinking $\overset\circ\to{\Bbb D}$ and $\Cal U$ if necessary, let $\tilde h: \overset\circ\to{\Bbb D}\times\Cal
U\rightarrow\Bbb C$ denote a local extension of $h$ to $\overset\circ\to{\Bbb D}\times\Cal U$. We use $w$ as our
coordinate on $\overset\circ\to{\Bbb D}$. Note that replacing $h(w, \bold z)$ by $h(w^2, \bold z)$  does not change the
statement of the theorem. Therefore, we can, and will, assume that $d_{(0, \bold x)}\tilde h$ vanishes on $\Bbb
C\times\{\bold 0\}$. 

\vskip .1in

Let $\tilde p:\overset\circ\to{\Bbb D}\times\Cal U\rightarrow\Cal U$ denote the projection, and let $p:=\tilde
p_{|_{\overset\circ\to{\Bbb D}\times X}}$. Let $\bold Q^\bullet:=p^*\Pdot[1]$; as $\Pdot$ is perverse, so is $\bold
Q^\bullet$. Let $Y:=(\overset\circ\to{\Bbb D}\times X)\cap V(h-b)$, and let $\widehat
w:Y\rightarrow\overset\circ\to{\Bbb D}$ denote the projection. Let $\bold R^\bullet:=\bold Q^\bullet_{|_{Y}}[-1]$. Let
$\hat f:Y\rightarrow\Bbb C$ be given by $\hat f(w, \bold z):= f(\bold z)$. As we already know that the theorem is true
for ordinary families of functions, we wish to apply it to the family of functions $\hat f_{\hat w}$ and the family of
sheaves $\bold R^\bullet_{\hat w}$; this would clearly prove the desired result.

\vskip .1in

Thus, we need to prove two things: that $\bold R^\bullet$ is perverse near $(0, \bold x)$, and that the family $\bold
R^\bullet_{\hat w}$ is continuous at $(0, \bold x)$.

\vskip .2in 

Let $\{S_\alpha\}$ be a Whitney stratification, with connected strata, of $X$ with respect to which $\Pdot$ is
constructible. Refining the stratification if necessary, assume that $V(t-b)$ is a union of strata. Let
$\operatorname{Ch}(\Pdot)=\sum m_\alpha\conc$. Clearly, $\bold Q^\bullet$ is constructible with respect to the Whitney
stratification $\{\overset\circ\to{\Bbb D}\times S_\alpha\}$, and the characteristic cycle of $\bold Q^\bullet$ in
$T^*(\overset\circ\to{\Bbb D}\times\Cal U)$ is given by $\operatorname{Ch}(\bold Q^\bullet) = -\sum m_\alpha
\Big[\overline{T^*_{\overset\circ\to{\Bbb D}\times S_\alpha}(\overset\circ\to{\Bbb D}\times\Cal U)}\Big]$. 

Note that, for all $(\bold z,
\eta)\in T^*\Cal U$, 
$(\bold z,
\eta)\in\con$ if and only if $(0,\bold z, \eta\circ d_{(0, \bold z)}p)\in \overline{T^*_{\overset\circ\to{\Bbb D}\times
S_\alpha}(\overset\circ\to{\Bbb D}\times\Cal U)}$. As we are assuming that $d_{(0, \bold x)}\tilde h$ vanishes on $\Bbb
C\times\{\bold 0\}$ and that $h_0=t$, we know that $d_{(0,\bold x)}\tilde h= d_\bold x\tilde t\circ d_{(0, \bold
z)}\tilde p$. Thus, $(\bold x, d_\bold x\tilde t)\in\con$ if and only if
$(0,\bold x, d_{(0,\bold x)}\tilde h)\in \overline{T^*_{\overset\circ\to{\Bbb D}\times S_\alpha}(\overset\circ\to{\Bbb
D}\times\Cal U)}$. Therefore, $(\bold x, d_\bold x\tilde t)\in |\operatorname{Ch}(\Pdot)|$ if and only if $(0,\bold x,
d_{(0,\bold x)}\tilde h)\in |\operatorname{Ch}(\bold Q^\bullet)|$. As we are assuming that the family $\Pdot_t$ is
continuous at $\bold x$, we may apply Lemma 5.3  to conclude that $(\bold x, d_\bold x\tilde t)\not\in
|\operatorname{Ch}(\Pdot)|$ and, hence, $(0,\bold x, d_{(0,\bold x)}\tilde h)\not\in |\operatorname{Ch}(\bold
Q^\bullet)|$. It follows that, for all $(w, \bold z)$ near $(0, \bold x)$, $(w,\bold z, d_{(w,\bold z)}\tilde h)\not\in
|\operatorname{Ch}(\bold Q^\bullet)|$ and that the family
$\bold Q^\bullet_{h}$ is continuous at $(0, \bold x)$; that is, there exists an open neighborhood, $\Omega\times\Cal W$,
of $(0, \bold x)$ in
$\overset\circ\to{\Bbb D}\times\Cal U$, in which $\phi_{h-b}[-1]\bold Q^\bullet=0$ and such that, if $(w, \bold z)\in
\Omega\times\Cal W$ and 
$m_\alpha\neq 0$, then
$(w,\bold z, d_{(w,\bold z)}\tilde h)\not\in \overline{T^*_{\overset\circ\to{\Bbb D}\times
S_\alpha}\big(\overset\circ\to{\Bbb D}\times\Cal U\big)}$. For the remainder of the proof, we assume that
$\overset\circ\to{\Bbb D}$ and $\Cal U$ have been rechosen to be small enough to use for $\Omega$ and $\Cal W$. 

As $\phi_{h-b}[-1]\bold Q^\bullet=0$, $\bold R^\bullet\cong\psi_{h-b}[-1]\bold Q^\bullet$ is a perverse sheaf on
$Y$. It remains for us to show that the family $\bold R^\bullet_{\hat w}$ is continuous  at $(0,
\bold x)$. 

\vskip .1in

Of course, we appeal to Lemma 5.3 again -- we need to show that $(0, \bold x, d_{(0,
\bold x)}w)\not\in |\operatorname{Ch}(\bold R^\bullet)|$. Now, $|\operatorname{Ch}(\bold R^\bullet)| =
|\operatorname{Ch}(\psi_{h-b}[-1]\bold Q^\bullet)|$, and we wish to use Theorem 2.3 to describe this characteristic
variety. If $(w, \bold z)\in \Omega\times\Cal W$ and 
$m_\alpha\neq 0$, then
$(w,\bold z, d_{(w,\bold z)}\tilde h)\not\in \overline{T^*_{\overset\circ\to{\Bbb D}\times
S_\alpha}\big(\overset\circ\to{\Bbb D}\times\Cal U\big)}$; thus, if $m_\alpha\neq 0$, then $h$ has no critical points
when restricted to $\overset\circ\to{\Bbb D}\times S_\alpha$, and, using the notation of 2.2 and 2.3, 
$$ T^*_{h-b, \bold Q^\bullet}\big(\overset\circ\to{\Bbb D}\times\Cal U\big) =
\sum_{\alpha}m_\alpha\left[\overline{T^*_{h_{|_{\overset\circ\to{\Bbb D}\times S_\alpha}}}\big(\overset\circ\to{\Bbb
D}\times\Cal U\big)}\right].
$$ Now, using Theorem 2.3, we find that 
$$ |\operatorname{Ch}(\bold R^\bullet)| = \big(V(h-b)\times \Bbb C^{n+2}\big)\cap
\bigcup_{m_\alpha\neq 0}\overline{T^*_{h_{|_{\overset\circ\to{\Bbb D}\times S_\alpha}}}\big(\overset\circ\to{\Bbb
D}\times\Cal U\big)}.
$$ We will be finished if we can show that, if $m_\alpha\neq 0$, then $(0, \bold x, d_{(0,
\bold x)}w)\not\in \overline{T^*_{h_{|_{\overset\circ\to{\Bbb D}\times S_\alpha}}}\big(\overset\circ\to{\Bbb D}\times\Cal
U\big)}$. 

Fix an $S_\alpha$ for which $m_\alpha\neq 0$. Suppose that $(0, \bold x, \eta)\in
\overline{T^*_{h_{|_{\overset\circ\to{\Bbb D}\times S_\alpha}}}\big(\overset\circ\to{\Bbb D}\times\Cal U\big)}$. Then,
there exists a sequence $(w_i, \bold z_i, \eta_i)\in T^*_{h_{|_{\overset\circ\to{\Bbb D}\times
S_\alpha}}}\big(\overset\circ\to{\Bbb D}\times\Cal U\big)$ such that $(w_i, \bold z_i,
\eta_i)\rightarrow(0, \bold x, \eta)$. Thus, $\eta_i\big((\Bbb C\times T_{\bold
z_i}S_\alpha)\cap\operatorname{ker}d_{(w_i,
\bold z_i)}\tilde h\big) = 0$. By taking a subsequence, if necessary, we may assume that $T_{\bold z_i}S_\alpha$
converges to some $\Cal T$ in the appropriate Grassmanian. Now, we know that $\operatorname{ker}d_{(w_i, \bold
z_i)}\tilde h\rightarrow \operatorname{ker}d_{(0, \bold x)}\tilde h = \Bbb C\times\operatorname{ker}d_\bold x\tilde t$.
As $(\bold x, d_\bold x\tilde t)\not\in\con$, $\Bbb C\times\operatorname{ker}d_\bold x\tilde t$ transversely intersects
$\Bbb C\times \Cal T$. Therefore, $(\Bbb C\times T_{\bold z_i}S_\alpha)\cap\operatorname{ker}d_{(w_i,
\bold z_i)}\tilde h\rightarrow (\Bbb C\times\Cal T)\cap (\Bbb C\times\operatorname{ker} d_\bold x\tilde t)$ , and so
$\Bbb C\times\{\bold 0\}\subseteq\operatorname{ker}\eta$. However, $\operatorname{ker} d_{(0,\bold x)} w =
\{0\}\times\Bbb C^{n+1}$, and we are finished.\qed

\vskip .4in

We would like to translate  Theorem 5.7 into a statement about Milnor fibres and the constant sheaf. First, though, it
will be convenient to prove a lemma.

\vskip .3in

\noindent{\bf Lemma 5.8}. {\it Let $\bold x\in X$, and let $b:=t(\bold x)$.  Suppose that ${\dm}_\bold x
\big(V(t-b)\cap\overline{\Sigma_{{}_\Bbb C}t}\big)\leqslant 0$. Fix an integer $k$. If $\widetilde H^k(F_{t, \bold x};\
\Bbb C) = 0$, then the family $\kp_t$ is continuous at
$\bold x$. In addition, if $\widetilde H^k(F_{t, \bold x};\ \Bbb C) = 0$ and $\widetilde H^{k-1}(F_{t, \bold x};\ \Bbb
C) = 0$, then
$\kp_b\cong {}^p\negmedspace H^0(\Bbb C_{X_b}^\bullet[k])$ near $\bold x$. }

\vskip .3in

\noindent{\it Proof}. By Remark 1.7, the assumption that ${\dm}_\bold x \big(V(t-b)\cap\overline{\Sigma_{{}_\Bbb
C}t}\big)\leqslant 0$ is equivalent to ${\dm}_\bold x\overline{\Sigma_{{}_\Bbb C}t}\leqslant 0$ and, by Theorem 3.6,
this is equivalent to ${\dm}_\bold x \overline{\Sigma_{{}_{{}^j\hskip -.01in\Pdot}}t}\leqslant 0$ for all $j$. Thus,
$\operatorname{supp}\phi_{t-b}[-1]\kp\subseteq\{\bold x\}$ near $\bold x$. We claim that the added
assumption that
$\widetilde H^k(F_{t, \bold x};\ \Bbb C) = 0$ implies that, in fact, $\phi_{t-b}[-1]\kp = 0$ near
$\bold x$.

For, near $\bold x$, $\operatorname{supp}\phi_{t-b}[-1]\Bbb C_X^\bullet[k+1]\subseteq\{\bold x\}$, and so 
$$
\phi_{t-b}[-1]\kp = \phi_{t-b}[-1]{}^p\negmedspace H^0(\Bbb C_X^\bullet[k+1])\cong {}^p\negmedspace
H^0(\phi_{t-b}[-1]\Bbb C_X^\bullet[k+1])\cong \bold H^0(\phi_{t-b}[-1]\Bbb C_X^\bullet[k+1]).
$$
 Near $\bold x$,  $\phi_{t-b}[-1]\Bbb C_X^\bullet[k+1]$ is supported at, at most, the point $\bold x$ and, hence,
$\phi_{t-b}[-1]\kp = 0$ provided that $\bold H^0(\phi_{t-b}[-1]\Bbb C_X^\bullet[k+1])_\bold x = 0$,
i.e., provided that $\widetilde H^k(F_{t, \bold x};\ \Bbb C) = 0$. This proves the first claim in the lemma.

\vskip .1in

Now, if the family $\kp_t$ is continuous at $\bold x$, then, near $\bold x$,
$$ \kp_b = \kp_{|_{V(t-b)}}[-1]\cong\psi_{t-b}[-1]{}^p\negmedspace H^0(\Bbb
C_X^\bullet[k+1])\cong {}^p\negmedspace H^0(\psi_{t-b}[-1]\Bbb C_X^\bullet[k+1]),
$$ and we claim that, if $\widetilde H^k(F_{t, \bold x};\ \Bbb C) = 0$ and $\widetilde H^{k-1}(F_{t, \bold x};\ \Bbb C)
= 0$, then there is an isomorphism (in the derived category) ${}^p\negmedspace H^0(\psi_{t-b}[-1]\Bbb
C_X^\bullet[k+1])\cong {}^p\negmedspace H^0(\Bbb C_{X_b}^\bullet[k])$.

To see this, consider the canonical distinguished triangle 
$$
\Bbb C_{X_b}^\bullet[k]\ \rightarrow\ \psi_{t-b}[-1]\Bbb C_X^\bullet[k+1]\ \rightarrow\ \phi_{t-b}[-1]\Bbb
C_X^\bullet[k+1]\  @>[1]>>\ \Bbb C_{X_b}^\bullet[k].
$$ A portion of the long exact sequence (in the category of perverse sheaves) resulting from applying perverse
cohomology is given by
$$ {}^p\negmedspace H^{-1}(\phi_{t-b}[-1]\Bbb C_X^\bullet[k+1])\rightarrow{}^p\negmedspace H^{0}(\Bbb
C_{X_b}^\bullet[k])\rightarrow{}^p\negmedspace H^0(\psi_{t-b}[-1]\Bbb C_X^\bullet[k+1])\rightarrow{}^p\negmedspace
H^0(\phi_{t-b}[-1]\Bbb C_X^\bullet[k+1]).
$$ We would be finished if we knew that the terms on both ends of the above were zero. However, since $\phi_{t-b}[-1]\Bbb
C_X^\bullet[k+1]$ has no support other than $\bold x$ (near $\bold x$), we proceed as we did above to show that 
${}^p\negmedspace H^{-1}(\phi_{t-b}[-1]\Bbb C_X^\bullet[k+1])$ and ${}^p\negmedspace H^{0}(\phi_{t-b}[-1]\Bbb
C_X^\bullet[k+1])$ are zero precisely when $\widetilde H^{k-1}(F_{t, \bold x};\ \Bbb C)$ and $\widetilde H^{k}(F_{t,
\bold x};\ \Bbb C)$ are zero.
\qed

\vskip .4in

\noindent{\bf Theorem 5.9}. {\it  Let $\bold x\in X$ and let $b:= t(\bold x)$.  Suppose that 
$\bold x\not\in\overline{\Sigma_{{}_\Bbb C}t}$, and that $\dm_\bold
x\Sigma_{{}_{\Bbb C}}(f_b)\leqslant 0$.

Then, 
there exists a neighborhood,
$\Cal W$, of $\bold x$ in $X$ such that, for all $a$ near
$b$, there are a finite number of points $\bold y\in\Cal W\cap V(t-a)$ for which $\widetilde H^*(F_{f_a, \bold y};\ \Bbb
C)\neq 0$; moreover, for all integers, $k$, $\tilde b_{k-1}(F_{f_a, \bold y})=\mu_{\bold y}(f_a; \kp_a)$, and
$$
\tilde b_{k-1}(F_{f_b, \bold x}) = \sum_{\bold y\in\Cal W\cap V(t-a)}\tilde b_{k-1}(F_{f_a, \bold y}),
$$ where $\widetilde H^*()$ and $\tilde b_*()$ are as in  Remark 3.4. }

\vskip .3in

\noindent{\it Proof}.  Let
$v:=f_b(\bold x)$. Fix an integer $k$. 

 By the lemma, the family $\kp_t$ is continuous at $\bold x$ and $\kp_b\cong
{}^p\negmedspace H^0(\Bbb C_{X_b}^\bullet[k])$ near $\bold x$. Thus, 
$$\phi_{f_b-v}[-1]\kp_b\cong \phi_{f_b-v}[-1]{}^p\negmedspace H^0(\Bbb C_{X_b}^\bullet[k])\cong
{}^p\negmedspace H^0(\phi_{f_b-v}[-1]\Bbb C_{X_b}^\bullet[k]).$$ We are assuming that $\dm_\bold
x\Sigma_{{}_{\Bbb C}}(f_b)\leqslant 0$; this is equivalent to: 
 $\operatorname{supp}\phi_{f_b-v}[-1]\Bbb C_{X_b}^\bullet[k]\subseteq\{\bold x\}$ near $\bold x$, it follows from the
above line and Theorem 3.6 that $\dm_\bold x\Sigma_{{}_{\kp_b}}f_b\leqslant 0$ and that 
$$\mu_\bold x(f_b; \kp_b) = \dm H^0(\phi_{f_b-v}[-1]\Bbb C_{X_b}^\bullet[k])_\bold x = \tilde
b_{k-1}(F_{f_b, \bold x}).\tag{$\ddagger$}$$

Applying Theorem 5.7, we find that there exists an open neighborhood $\Cal W^\prime$ of $\bold x$ in $\Cal U$ such
that,  for all $\bold y\in\Cal W^\prime$, if $a:=t(\bold y)$, then  $(*)$ $\dm_\bold y\Sigma_{{}_{\kp_a}}f_a\leqslant 0$, and, for fixed $a$ close to $b$, there are a finite number of points $\bold y\in\Cal
W^\prime\cap V(t-a)$ such that  $\mu_\bold y(f_a; \kp_a)\neq 0$ and
$$
\mu_\bold x(f_b; \kp_b) = \sum_{\bold y\in\Cal W\cap V(t-a)}\mu_\bold y(f_a; \kp_a).\tag{$\dagger$}
$$

Now, using the above argument for all $k$ with $0\leqslant k\leqslant \dm X-1$  and intersecting the resulting
$\Cal W^\prime$-neighborhoods, we obtain an open neighborhood $\Cal W$ of $\bold x$ such that $(*)$ and $(\dagger)$ hold
for all such $k$. We claim that, if $a$ is close to
$b$, then $\Cal W\cap\overline{\Sigma_{\Bbb C}f_a}$ consists of isolated points, i.e., the points $\bold y\in\Cal W\cap
V(t-a)$ for which $\widetilde H^*(F_{f_a, \bold y};\ \Bbb C)\neq 0$ are isolated.

If $a=b$, then there is nothing to show. So, assume that $a\neq b$, and assume that we are working in $\Cal W$
throughout. By Remark 1.7,
$t$ satisfies the hypotheses of Lemma 5.8 at
$t=a$; hence, for all $k$, not only is $\kp_t$ continuous at $t=a$, but we also know that $\kp_a\cong {}^p\negmedspace H^0(\Bbb C_{X_a}^\bullet[k])$. By Theorem 3.6, $\overline{\Sigma_{\Bbb
C}f_a}=\bigcup
\overline{\Sigma_{\kp_a}f_a}$, where the union is over $k$ where $0\leqslant k\leqslant \dm X_a$. As
$\dm X_a\leqslant\dm X -1$, the claim follows from $(*)$ and the definition of $\Cal W$. 

\vskip .1in

Now that we know that $\kp_t$ is continuous at $t=a$ and that $\Cal W\cap\overline{\Sigma_{\Bbb
C}f_a}$ consists of isolated points, we may use the argument that produced $(\ddagger)$ to conclude that $\mu_\bold
y(f_a; \kp_a)  = \tilde b_{k-1}(F_{f_a, \bold y})$. The theorem follows from this, $(\ddagger)$,
$(*)$, and
$(\dagger)$.
\qed

\vskip .4in

We want to prove a result which generalizes that of L\^e and Saito [{\bf L-S}]. We need to make the assumption that the
Milnor number is constant along a curve that is embedded in $X$. Hence, it will be convenient to use a {\it local
section} of
$t:X\rightarrow \Bbb C$ at a point $\bold x\in X$; that is, an analytic function $\bold r$ from an open neighborhood,
$\Cal V$, of $t(\bold x)$ in $\Bbb C$ into $X$ such that $\bold r(t(\bold x)) =
\bold x$ and $t\circ\bold r$ equals the inclusion morphism of $\Cal V$ into $\Bbb C$. Note that existence of such a
local section implies that $\bold x\not\in\Sigma_{\operatorname{alg}}t$; in particular, $V(\tilde t-\tilde t(\bold x))$
is smooth at $\bold x$.

\vskip .3in

\noindent{\bf Theorem 5.10}. {\it  Suppose that the family $\Pdot_t$ is continuous at
$\bold x\in X$. Let $b:= t(\bold x)$, and let $v:=f_b(\bold x)$. Let $\bold r:\Cal V\rightarrow X$ be a local section of
$t$ at $\bold x$, and let $C:=\operatorname{im}\bold r$. Assume that $C\subseteq V(f-v)$, that 
$\dm_\bold x\Sigma_{{}_{\Pdot_b}}f_b\leqslant 0$, and that, for all $a$ close to $b$, the Milnor number $\mu_{\bold
r(a)}(f_a; \Pdot_a)$ is non-zero and is independent of $a$; denote this common value by $\mu$.

\vskip .1in

Then, $C$ is smooth at $\bold x$, $V(\tilde t-b)$ transversely intersects $C$ in $\Cal U$ at $\bold x$ , and there
exists a neighborhood, $\Cal W$, of $\bold x$ in
$X$ such that $\Cal W\cap\overline{\Sigma_{\Pdot}f}\subseteq C$ and $\big(\phi_{f-v}[-1]\Pdot\big)_{|_{\Cal W\cap
C}}\cong \big(\Bbb C^\mu_{\Cal W\cap C}[1]\big)^{\bullet}$. In particular, if we let $\hat t$ denote the restriction of
$t$ to $V(f-v)$, then the family $\big(\phi_{f-v}[-1]\Pdot\big)_{\hat t}$ is continuous at $\bold x$.

\vskip .1in

If, in addition to the other hypotheses, we assume that $\bold x\in\Sigma_{\operatorname{alg}}f$, then the two families
$\big(\psi_{f-v}[-1]\Pdot\big)_{\hat t}$ and $\big(\Pdot_{|_{V(f-v)}}[-1]\big)_{\hat t}$ are continuous at $\bold x$.
(Though
$\Pdot_{|_{V(f-v)}}[-1]$ need not be perverse.) }

\vskip .3in

\noindent{\it Proof}.  Let us first prove that the last statement of the theorem follows easily from the first portion
of the theorem. So, assume that  $\phi_{\hat t-b}[-1]\phi_{f-v}[-1]\Pdot = 0$ near $\bold x$. Therefore, working near
$\bold x$, we have that $\phi_{\hat t-b}[-1]\big(\Pdot_{|_{V(f-v)}}[-1]\big)\cong \phi_{\hat
t-b}[-1]\psi_{f-v}[-1]\Pdot$, and we need to show that this is the zero-sheaf. By Lemma 5.3, what we need to show is
that $(\bold x, d_\bold x\tilde t)\not\in|\operatorname{Ch}(\psi_{f-v}[-1]\Pdot)| =
|\operatorname{Ch}(\psi_{f-v}\Pdot)|$. As we are assuming that $\bold x\in\Sigma_{\operatorname{alg}}f$, we may apply
Corollary 4.5 to find that it suffices to show that $(\bold x, d_\bold x\tilde
t)\not\in|\operatorname{Ch}(\phi_{f-v}\Pdot)|= |\operatorname{Ch}(\phi_{f-v}[-1]\Pdot)|$. By 5.3, this is equivalent to
$\phi_{\hat t-b}[-1]\phi_{f-v}[-1]\Pdot = 0$ near $\bold x$, which we already know to be true. This proves the last
statement of the theorem.

\vskip .2in

Before proceeding with the remainder of the proof, we wish to make some simplifying assumptions. As $\bold
x\not\in\Sigma_{\operatorname{alg}}t$, we may certainly perform an analytic change of coordinates in
$\Cal U$ to reduce ourselves to the case where $t$ is simply the restriction to $X$ of a linear form $\tilde t$.
Moreover, it is notational convenient to assume, without loss of generality, that $\bold x=\bold 0$ and that $b$ and $v$
are both zero.

\vskip .1in

Let $\{S_\alpha\}$ be a Whitney stratification of $X$ with connected strata with respect to which $\Pdot$ is
constructible and such that
$V(t)$ and $V(f)$ are each unions of strata. Suppose that $\operatorname{Ch}(\Pdot)$ is given by $\sum_\alpha
m_\alpha\conc$. 

\vskip .1in

Let $\widetilde C:= \{(\bold r(a), d_{\bold r(a)}\tilde f)\ |\ a\in\Cal V\}$; the projection, $\rho$, onto the first
component induces an isomorphism from $\widetilde C$ to $C$. By Lemma 5.6, the assumption that the Milnor number,
$\mu_{\bold r(a)}(f_a; \kp_a)$, is independent of $a$ is equivalent to: 

\vskip .1in

\noindent $(\dagger)$ there exists an open neighborhood $\widetilde{\Cal W}$ of $(\bold 0, d_\bold 0\tilde f)$ in
$T^*\Cal U$ in which
$\widetilde C$ equals 
$$\imdf\ \cap\ \bigcup_{m_\alpha\neq 0}\overline{T^*_{t_{|_{S_\alpha}}}\Cal U}$$ and $\widetilde C$
is a smooth curve at
$(\bold 0,d_\bold 0\tilde f)$ such that $(\bold 0, d_\bold 0\tilde f)\not\in\Sigma(t\circ\rho_{|_{\widetilde C}})$. 
\vskip .2in

\noindent It follows immediately that
$C$ is smooth at $\bold 0$ and $\bold 0\not\in\Sigma(t_{|_{C}})$. We need to show that $(\dagger)$ implies that  $\Cal
W\cap\overline{\Sigma_{\Pdot}f}\subseteq C$ and
$\big(\phi_{f}[-1]\Pdot\big)_{|_{\Cal W\cap C}}\cong \big(\Bbb C^\mu_{\Cal W\cap C}[1]\big)^{\bullet}$, where $\Cal
W:=\rho(\widetilde{\Cal W})$.

\vskip .1in

As $\overline{T^*_{S_\alpha}\Cal U}\subseteq \overline{T^*_{t_{|_{S_\alpha}}}\Cal U}$, we have that
$|\operatorname{Ch}(\Pdot)|\subseteq
\bigcup_{m_\alpha\neq 0}\overline{T^*_{t_{|_{S_\alpha}}}\Cal U}$ and, thus, $\imdf\cap
|\operatorname{Ch}(\Pdot)|\subseteq\widetilde C$ inside $\widetilde{\Cal W}$. It follows from Theorem 3.2 that 
$\Cal W\cap\overline{\Sigma_{\Pdot}f}\subseteq C$.

\vskip .1in

It remains for us to show that $\big(\phi_{f}[-1]\Pdot\big)_{|_{\Cal W\cap C}}\cong \big(\Bbb C^\mu_{\Cal W\cap
C}[1]\big)^{\bullet}$. As $\phi_f[-1]\Pdot$ is perverse and we have just shown that the support of $\phi_{f}[-1]\Pdot$,
near $\bold 0$, is a smooth curve, it follows from the work of MacPherson and Vilonen in [{\bf M-V}] that what we need
to show is that, for a generic linear form $L$, $\bold Q^\bullet:=\phi_{L}[-1]\phi_{f}[-1]\Pdot = 0$ near $\bold 0$.  By
definition of the characteristic cycle (and since $\bold 0$ is an isolated point in the support of $\bold Q^\bullet$),
this is the same as showing that the coefficient of
$T^*_{\{\bold 0\}}\Cal U$ in 
$\operatorname{Ch}(\phi_{f}[-1]\Pdot)$ equals zero. To show this, we will appeal to Theorem 2.4 and use the notation
from there.

We need to show that $m_\bold 0(\phi_{f}[-1]\Pdot)=0$. By 2.4, if suffices to show that $m_\bold 0(\Pdot)=0$ and
$\Gamma^1_{f,L}(S_\alpha)=\emptyset$ near $\bold 0$, for all $S_\alpha$ which are not contained in $V(f)$ and for which
$m_\alpha\neq 0$ (where
$L$ still denotes a generic linear form). As $\Pdot_t$ is continuous at $\bold 0$, Lemma 5.3 tells us that $m_\bold
0(\Pdot)=0$. Now, near
$\bold 0$, if
$\bold y\in \Gamma^1_{f,L}(S_\alpha)-\{\bold 0\}$, then $(\bold y, d_\bold y\tilde f)\in T^*_{L_{|_{S_\alpha}}}\Cal U$.
If we knew that, near $(\bold 0, d_\bold 0\tilde f)$, $\widetilde C$ equals 
$\imdf\ \cap\ \bigcup_{m_\alpha\neq 0}\overline{T^*_{L_{|_{S_\alpha}}}\Cal U}$, then we would be
finished -- for
$C$ is contained in $V(f)$ while $S_\alpha$ is not; hence, $\Gamma^1_{f,L}(S_\alpha)$ would have to be empty near $\bold
0$. 

Looking back at $(\dagger)$, we see that what we still need to show is that if $\widetilde C$ equals 
$\imdf\ \cap\ \bigcup_{m_\alpha\neq 0}\overline{T^*_{t_{|_{S_\alpha}}}\Cal U}$ near $(\bold 0,
d_\bold 0\tilde f)$, then the same statement holds with $t$ replaced by a generic linear form $L$. We accomplish this by
perturbing $t$ until it is generic, and by then showing that this perturbed $t$ satisfies the hypotheses of the theorem.

As $C$ is smooth and transversely intersected by $V(\tilde t)$ at $\bold 0$, by performing an analytic change of
coordinates, we may assume that $\tilde t=z_0$, that $C$ is the $z_0$-axis, and that $r(a)= (a, \bold 0)$. Since the set
of linear forms for which 2.4 holds is generic, there exists an open disk, $\overset\circ\to{\Bbb D}$, around the origin
in $\Bbb C$ and an analytic family $\tilde h:(\overset\circ\to{\Bbb D}\times\Cal U, \overset\circ\to{\Bbb
D}\times\{\bold 0\})\rightarrow(\Bbb C, 0)$ such that $\tilde h_0(\bold z):= \tilde h(0, \bold z) = \tilde t(\bold z)$
and such that, for all small non-zero $c$, $\tilde h_c(\bold z):=\tilde h(c, \bold z)$ is a linear form for which
Theorem 2.4 holds. Let $h:=\tilde h_{|_{\overset\circ\to{\Bbb D}\times X}}$.

As the family $\Pdot_t$ is continuous at $\bold 0$, Lemma 5.5 tells us that $\Pdot_{h_c}$ is continuous at $\bold 0$ for
all small $c$. As we are now considering these two different families with the same underlying sheaf, the expression
$\Pdot_a$ for a fixed value of
$a$ is ambiguous, and we need to adopt some new notation. We continue to let $\Pdot_a:= \Pdot_{|_{V(t-a)}}[-1]$ and
$f_a:=f_{|_{V(t-a)}}$, and let
${}_c\Pdot_a:=\Pdot_{|_{V(h_c-a)}}[-1]$ and ${}_cf_a:=f_{|_{V(h_c-a)}}$.  

Since $V(\tilde h_0)=V(z_0)$  transversely intersects $C$ at $\bold 0$ in $\Cal U$, for all small $c$, $V(h_c)$
transversely intersects $C$ at
$\bold 0$ in $\Cal U$. Hence, for all small $c$, there exists a local section $\bold r_c(a)$ for $h_c$ at $\bold 0$ such
that
$\operatorname{im}\bold r_c\subseteq C$. 

\vskip .1in

We claim that, for all small $c$:

\vskip .1in

\noindent i)\hskip .2in $\dm_\bold 0\Sigma_{{}_c\Pdot_0}({}_cf_0)\leqslant 0$ and $\mu_\bold 0({}_cf_0;\
{}_c\Pdot_0)\leqslant \mu_\bold 0({}_0f_0;\ {}_0\Pdot_0) = \mu_\bold 0(f_0;\ \Pdot_0)$;

\vskip .1in

\noindent ii)\hskip .16in for all small $a$, $\dm_{\bold r_c(a)}\Sigma_{{}_c\Pdot_a}({}_cf_a)\leqslant 0$ and
$\mu_{\bold r_c(a)}({}_cf_a;\ {}_c\Pdot_a)\leqslant \mu_\bold 0({}_0f_0;\ {}_0\Pdot_0)$; and

\vskip .1in

\noindent iii)\hskip .12in for all small $a\neq 0$, $\mu_{\bold r_c(a)}({}_cf_a;\ {}_c\Pdot_a) = \mu_{\bold
r_c(a)}(f_{z_0(\bold r_c(a))};\
\Pdot_{z_0(\bold r_c(a))})$. 

\vskip .2in

\noindent Note that proving i), ii), and iii) would complete the proof of the theorem, for they imply that the
hypotheses of the theorem hold with $t$ replaced by $h_c$ for all small $c$. To be precise, we would know that
$\Pdot_{h_c}$ is continuous at $\bold 0$, $\dm_\bold 0\Sigma_{{}_c\Pdot_0}({}_cf_0)\leqslant 0$, and, for all small $a$,
$\mu_{\bold r_c(a)}({}_cf_a;\ {}_c\Pdot_a)= \mu_{\bold 0}({}_cf_0;\ {}_c\Pdot_0)$; this last equality follows from i),
ii), and iii), since, for all small $a\neq 0$, we would have 
$$\mu=\mu_{\bold r_c(a)}(f_{z_0(\bold r_c(a))};\
\Pdot_{z_0(\bold r_c(a))})= \mu_{\bold r_c(a)}({}_cf_a;\ {}_c\Pdot_a)\leqslant \mu_{\bold 0}({}_cf_0;\
{}_c\Pdot_0)\leqslant \mu_{\bold 0}(f_0;\ \Pdot_0)=\mu.$$

However, i), ii) and iii) are easy to prove. i) and ii) follow immediately from Theorem 5.7, and iii) follows simply
from the fact that, for all small $a\neq0$, $V(z_0-z_0(\bold r_c(a)))$ and $V(\tilde h_c-\tilde h_c(\bold r_c(a)))$ are
smooth and transversely intersect all strata of any analytic stratification of $X$ in a neighborhood of $(0, \bold 0)$.
This concludes the proof.\qed 

\vskip .4in

\noindent{\bf Corollary 5.11}. {\it  Suppose that the family $\Pdot_t$ is continuous at
$\bold x\in X$. Let $b:= t(\bold x)$, and let $v:=f_b(\bold x)$. Let $\bold r:\Cal V\rightarrow X$ be a local section of
 $t$ at $\bold x$, and let $C:=\operatorname{im}\bold r$. Assume that $C\subseteq V(f-v)$, that 
$\dm_\bold x\Sigma_{{}_{\Pdot_b}}f_b\leqslant 0$, and that, for all $a$ close to $b$, the Milnor number $\mu_{\bold
r(a)}(f_a; \Pdot_a)$ is non-zero and is independent of $a$. Let $\operatorname{Ch}(\Pdot) =\sum_\alpha m_\alpha\conc$,
where
$\{S_\alpha\}$ is a collection of connected analytic submanifolds of $\Cal U$.

\vskip .1in

Then, $C$ is smooth at $\bold x$, and there exists a neighborhood, $\Cal W$, of $\bold x$ in
$X$ such that, for all $S_\alpha$ for which $S_\alpha\not\subseteq V(f-v)$ and $m_\alpha\neq 0$:

 $\Cal W\cap\Sigma\big(f_{|_{(\overline{S_\alpha})_{\operatorname{reg}}}}\big)\subseteq C$ and, if $\bold
x\in\Sigma_{{}_{\operatorname{Nash}}}(f_{|_{{\overline S}_\alpha}})$, then the pair
$\Big(\big(\overline{S_\alpha}\big)_{\operatorname{reg}},\
C\Big)$ satisfies Thom's
$a_f$ condition at $\bold x$. }
\vskip .3in

\noindent{\it Proof}. One applies Theorem 5.10. The fact that 
$\Cal W\cap\Sigma\big(f_{|_{(\overline{S_\alpha})_{\operatorname{reg}}}}\big)\subseteq C$, for all $S_\alpha$ for which
$m_\alpha\neq 0$  follows from Theorem 3.2, since $\Cal W\cap\overline{\Sigma_{\Pdot}f}\subseteq C$. The remainder of
the corollary follows by applying Corollary 4.4, where one uses $C$ for the submanifold
$N$.\qed

\vskip .4in

Just as we used perverse cohomology to translate Theorem 5.7 into a statement about the constant sheaf in Theorem 5.9,
we can use perverse cohomology to translate Corollary 5.11. We will use the notation and results from Proposition 3.3
and Remark 3.4.

\vskip .3in

\noindent{\bf Corollary 5.12}. {\it  Let $b:= t(\bold x)$, and let $v:=
f_b(\bold x)$.  Suppose that $\bold x\not\in\overline{\Sigma_{{}_\Bbb C}t}$. Suppose, further, that, $\dm_\bold
x\Sigma_{{}_{\Bbb C}}(f_b)\leqslant 0$.

Let $\bold r:\Cal V\rightarrow X$ be a local section of $t$ at $\bold x$, and let $C:=\operatorname{im}\bold r$. Assume
that $C\subseteq V(f-v)$.

\vskip .1in

Let $S_\alpha$ be a visible stratum of $X$ of dimension $d_\alpha$, not contained in $V(f-v)$, and let $j$ be an integer such that
$\tilde b_{j-1}(\Bbb L_\alpha)\neq 0$. Let $Y:= \overline{S_\alpha}$ and let $k:=d_\alpha+j-1$. In particular, $Y$
could be any irreducible component of $X$, $j$ could be zero, and $k$ would be $(\dm Y)-1$.

\vskip .1in

Suppose that the reduced Betti number $\tilde
b_{k-1}(F_{f_a, \bold r(a)})$ is independent of $a$ for all small $a$, and that either

\vskip .1in

a)\hskip .2in  $\bold x\in \Sigma_{{}_{\operatorname{Nash}}}(f_{|_Y})$; or that

\vskip .1in

b)\hskip .2in $\bold x\not\in \Sigma_{\operatorname{cnr}}(f_{|_Y})$, $C$ is smooth at $\bold x$, and
$(Y_{\operatorname{reg}}, C)$ satisfies Whitney's condition a) at
$\bold x$.

\vskip .2in

\noindent Then, $C$ is smooth at $\bold x$, and the pair $(Y_{\operatorname{reg}},\ C)$ satisfies the $a_f$ condition at $\bold x$.

\vskip .2in

Moreover, in case a), $\tilde
b_{k-1}(F_{f_a, \bold r(a)})\neq 0$, $C$ is transversely intersected
by $V(\tilde t-b)$ at $\bold x$,  and
$\Sigma(f_{|_{Y_{\operatorname{reg}}}})\subseteq C$ near $\bold x$. 
\vskip .2in

In addition, if $\bold x\in \Sigma_{\operatorname{alg}}f$ and, for all small $a$ and for all $i$, $\tilde
b_{i}(F_{f_a, \bold r(a)})$ is independent of $a$, then
$\bold x\not\in\overline{\Sigma_{{}_\Bbb C}(t_{|_{V(f-v)}})}$.

}
\vskip .3in

\noindent{\it Proof}.  We will dispose of case b) first. Suppose that $\bold x\not\in
\Sigma_{\operatorname{cnr}}(f_{|_Y})$, $C$ is smooth at $\bold x$, and $(Y_{\operatorname{reg}}, C)$ satisfies
Whitney's condition a) at $\bold x$. Let $\overset\circ\to
Y:=Y_{\operatorname{reg}}-\Sigma(f_{|_{Y_{\operatorname{reg}}}})$. 

Suppose that we have an analytic path $(\bold x(t), \eta_t)\in \overline{T^*_{f_{|_{\overset\circ\to Y}}}\Cal U}$,
where
$(\bold x(0),
\eta_0)= (\bold x, \eta)$ and, for $t\neq 0$, $(\bold x(t), \eta_t)\in T^*_{f_{|_{\overset\circ\to Y}}}\Cal U$. We wish
to show that  $(\bold x, \eta)\in T^*_{{}_C}\Cal U$. 

For $t\neq 0$, 
$\bold x(t)\in \overset\circ\to Y$, and thus 
$\eta_t$ can be written uniquely as $\eta_t=\omega_t+\lambda_t d_{\bold x_t}\tilde f$, where $\omega_t\in
T^*_{\overset\circ\to Y}\Cal U$ and $\lambda_t\in\Bbb C$. As we saw in Theorem 4.2, this implies that either
$|\lambda_t|\rightarrow\infty$ or that $\lambda_t\rightarrow\lambda_0$, for some $\lambda_0\in\Bbb C$. If
$|\lambda_t|\rightarrow\infty$, then $\frac{\eta_t}{\lambda_t}\rightarrow 0$ and, therefore,
$-\frac{\omega_t}{\lambda_t}\rightarrow d_{\bold x}\tilde f$; however, this implies that $\bold x\in
\Sigma_{\operatorname{cnr}}(f_{|_Y})$, contrary to our assumption. Thus, we must have that
$\lambda_t\rightarrow\lambda_0$.

It follows at once that $\omega_t$ converges to some $\omega_0$. By Whitney's condition a), $(\bold x, \omega_0)\in
T^*_{{}_C}\Cal U$. As $C\subseteq V(f-v)$, $(\bold x, d_{\bold x}\tilde f)\in T^*_{{}_C}\Cal U$. Hence, $(\bold x,
\eta)\in T^*_{{}_C}\Cal U$ and we have finished with case b).

\vskip .2in

We must now prove the results in case a). The main step is to prove that $\tilde
b_{k-1}(F_{f_b, \bold x})\neq 0$. 

\vskip .1in

We may refine our stratification, if necessary, so that $V(t-b)$ is a union of strata. By the first part of Theorem 5.9,
$\tilde b_{k-1}(F_{f_b,
\bold x})=\mu_\bold x(f_b; \kp_b)$. Hence, by Lemma 5.6.iii, $\tilde b_{k-1}(F_{f_b, \bold x})$
would be unequal to zero if we knew, for some $S_\beta$ for which $m_\beta\big(\kp\big)\neq
0$, that $(\bold x, d_{\bold x}\tilde f)\in\overline{T^*_{t_{|_{S_\beta}}}\Cal U}$. However, our fixed $S_\alpha$ is
such a stratum, for $b_{k+1-d_\alpha}(N_\alpha, \Bbb L_\alpha)\neq 0$ and, since $\bold x\in
\Sigma_{{}_{\operatorname{Nash}}}(f_{|_Y})$, $\bold x\in
\Sigma_{\operatorname{cnr}}(f_{|_Y})$ and so  $(\bold x, d_{\bold x}\tilde f)\in\overline{T^*_{S_\alpha}\Cal
U}\subseteq \overline{T^*_{t_{|_{S_\alpha}}}\Cal U}$.

\vskip .1in

Now, applying the first part of 5.9 again, we have that $\mu_{\bold r(a)}(f_a;\ \kp_a) = \tilde
b_{k-1}(F_{f_a, \bold r(a)})$ for all small $a$. The conclusions in case a) follow from Corollary 5.11.

\vskip .2in

We must still demonstrate the last statement of corollary.

\vskip .1in

Suppose that if $\tilde b_{i}(F_{f_a, \bold r(a)})$ is independent of $a$ for all small $a$ and for all $i$.  Let $\hat
t$ denote the restriction of $t$ to $V(f-v)$. We will work in a small neighborhood of $\bold x$. Applying the last two
sentences of Theorem 5.10, we find that $\phi_{\hat t-b}[-1]\phi_{f-v}[-1]{}^i\hskip -.01in\Pdot=0$ and $\phi_{\hat
t-b}[-1]\psi_{f-v}[-1]{}^i\hskip -.01in\Pdot=0$ for all $i$. Commuting nearby and vanishing cycles with perverse
cohomology, we find that 
$$ {}^p\negmedspace H^0\big(\phi_{\hat t-b}[-1]\phi_{f-v}[-1]\Bbb C_X^\bullet[i+1]\big) = 0\hskip .2in \text{ and
}\hskip .2in {}^p\negmedspace H^0\big(\phi_{\hat t-b}[-1]\psi_{f-v}[-1]\Bbb C_X^\bullet[i+1]\big) = 0,
$$ for all $i$. Therefore, $\phi_{\hat t-b}[-1]\phi_{f-v}[-1]\Bbb C_X^\bullet = 0$ and $\phi_{\hat
t-b}[-1]\psi_{f-v}[-1]\Bbb C_X^\bullet = 0$. It follows from the existence of the distinguished triangle (relating
nearby cycles, vanishing cycles, and restriction to the hypersurface) that
$\phi_{\hat t-b}[-1]\Bbb C_{V(f-v)}^\bullet[-1] = 0$. This proves the last statement of the corollary.\qed

\vskip .4in

\noindent{\it Remark 5.13}. If $X$ is a connected l.c.i., then each $\Bbb L_\alpha$ has (possibly) non-zero cohomology
concentrated in middle degree. Hence, for each visible $S_\alpha$, $\tilde b_{j-1}(\Bbb L_\alpha)\neq 0$ only when
$j={\operatorname{codim}}_{{}_X} S_\alpha$; this corresponds to $k=(\dm X)-1$. Therefore, the degree $(\dm X)-2$ reduced
Betti number of $F_{f_a, \bold r(a)}$ controls the $a_f$ condition between {\bf all} visible strata and $C$.

\vskip .4in

\noindent{\bf Corollary 5.14}. {\it   Let $W$ be an analytic subset of an open subset of $\Bbb C^n$. Let
$Z$ be a $d$-dimensional irreducible component of $W$. Let $X:=\overset\circ\to{\Bbb D}\times W$ be the product of an open disk about the
origin with $W$, and let $Y:=\overset\circ\to{\Bbb D}\times Z$. Let $f:(X, \overset\circ\to{\Bbb D}\times\{\bold 0\})\rightarrow (\Bbb C, 0)$
be an analytic function, such that $f_{|_Y}\not\equiv 0$, and let $f_t(\bold z):= f(t, \bold z)$. 

Suppose that $\bold 0$ is an isolated point of $\Sigma_{{}_{\Bbb
C}}(f_0)$, and that the reduced Betti number $\tilde b_{d-1}(F_{f_a, (a,\bold 0)})$ is independent of $a$ for all
small $a$.

If either \hskip .1in a) $\bold 0\in\Sigma_{{}_{\operatorname{Nash}}}(f_{|_Y})$ \ or \hskip .1in b) $\bold
0\not\in\Sigma_{\operatorname{cnr}}(f_{|_Y})$, then the pair $(Y_{\operatorname{reg}},\
\overset\circ\to{\Bbb D}\times\{\bold 0\})$ satisfies Thom's
$a_f$ condition at $\bold 0$. 

\vskip .2in

Moreover, in case a), $\tilde b_{d-1}(F_{f_a, (a,\bold 0)})\neq 0$ and, near
$\bold 0$, 
$\Sigma(f_{|_{Y_{\operatorname{reg}}}})\subseteq \overset\circ\to{\Bbb D}\times\{\bold 0\}$.

}

\vskip .5in

\noindent{\it Remark 5.15}. A question naturally arises: how effective is the criterion appearing in Corollary 5.14 that $\tilde
b_{d-1}(F_{f_a, (a,\bold 0)})$ is independent of $a$? 

By Proposition 3.10, if $\{R_\beta\}$ is a Whitney stratification of $W$, then (using the notation from 3.10)
$$\tilde b_{d-1}(F_{f_a, (a,\bold 0)})=$$
$$\tilde b_{d-1}(\Bbb L_{\{\bold 0\}})\ \ +\ \sum\Sb R_\beta\text{
visible}\\ \dm R_\beta\geqslant1\endSb\tilde b_{d-1-d_\beta}(\Bbb L_\beta)
\left(\big(\Gamma^1_{{}_{f_a,
L}}(R_\beta)\cdot V(f_a)\big)_\bold 0 - 
\big(\Gamma^1_{{}_{f_a, L}}(R_\beta)\cdot V(L)\big)_\bold 0
\right) ,
$$
where $\Bbb L_{\{\bold 0\}}$ denotes the complex link of the origin. As the Betti numbers do not vary with $a$, $\tilde b_{d-1}(F_{f_a,
(a,\bold 0)})$ will be independent of $a$ provided that $\big(\Gamma^1_{{}_{f_a,
L}}(R_\beta)\cdot V(f_a)\big)_\bold 0 - 
\big(\Gamma^1_{{}_{f_a, L}}(R_\beta)\cdot V(L)\big)_\bold 0$ is independent of $a$ for all visible strata, $R_\beta$, of dimension at least
one.

This condition is certainly very manageable to check if the dimension of the singular set of $X$ at the origin is zero or one.

\vskip .5in

The final statement of Corollary 5.12 has as its conclusion that the constant sheaf on
$X\cap V(f-v)$, parametrized by the restriction of  $t$, is continuous at $\bold x$; this is  useful for inductive
arguments, since the hypothesis on the ambient space in Corollary 5.12 is that the constant sheaf, parametrized by $t$,
should be continuous at $\bold x$. For instance, we can prove the following corollary.

\vskip .4in

\noindent{\bf Corollary 5.16}. {\it  Suppose that $f^1, \dots, f^k$  are analytic functions from $\Cal U$ into $\Bbb C$
which define a sequence of local complete intersections at the origin, i.e., are such that, for all $i$ with $1\leqslant
i\leqslant k$, the space
$X^{n+1-i}:=V(f^1, \dots, f^i)$ is a local complete intersection of dimension $n+1-i$ at the origin.  If, for all $i$,
$X^{n+1-i}_t$ has an isolated singularity at the origin  and the restrictions $f^{i+1}_t:X^{n+1-i}_t\rightarrow\Bbb C$
are such that
$\dm_\bold 0\Sigma_{\operatorname{can}}f^{i+1}_t\leqslant 0$ and have Milnor numbers (in the sense of {\rm [{\bf Lo}]})
which are independent of
$t$, then $\Sigma\big(f_{|_{X^{n+1-(k-1)}_{\operatorname{reg}}}}\big)\subseteq 
\Bbb C\times\{\bold 0\}$ and the pair $\big(X^{n+1-(k-1)}_{\operatorname{reg}},
\Bbb C\times\{\bold 0\}\big)$ satisfies the $a_{f^k}$ condition at the origin. }

\vskip .3in

\noindent{\it Proof}. Recall that $\Bbb C^\bullet_{X}[\dm X]$ is a perverse sheaf if $X$ is a local complete
intersection. The ``ordinary'' Milnor number of $f^{i+1}_t$ at the origin is equal to $\mu_\bold 0(f^{i+1}_t;\ \Bbb
C^\bullet_{X^{n+1-i}_t}[n-i])$.  Hence, using Proposition 3.10.ii, this Milnor number is equal to the degree $n-i-1$
(the ``middle'' degree)  reduced Betti number of the Milnor fibre of $f^{i+1}_t$ at the origin -- the only possible
non-zero reduced Betti number. Now, use Corollary 5.12 and induct; the inductive requirement on the Milnor fibre of
$z_0$ follows from the last statement of the corollary.\qed

\vskip .4in

\noindent{\it Remark 5.17}.  In [{\bf G-K}], Gaffney and Kleiman deal with families of local complete intersections as
above.  In this setting, they obtain the result of Corollary 5.16 using multiplicities of modules.

\vskip .4in

\noindent\S6. {\bf Concluding remarks}.

\vskip .2in

We hope to have convinced the reader that the correct notion of ``the critical locus'' of a function, $f$, on a singular
space  is given by
$\Sigma_{\Bbb C}f$.

\vskip .1in

We also hope to have convinced the reader of (at least) three other things: that the vanishing cycles control Thom's
$a_f$ condition (as demonstrated in Corollary 4.4, Corollary 5.11, and  Corollary 5.12), that the correct setting to be
in to generalize many classical results is where one uses arbitrary perverse sheaves as coefficients, and that perverse
cohomology is an amazing tool for turning statements about perverse sheaves into statements about the constant sheaf.

\vskip .1in

While a great deal of material concerning local complete intersections appears in the singularities literature, it is
not so easy to find results that apply to arbitrary analytic spaces. As we remarked earlier, from our point of view,
what is special about l.c.i.'s is that the shifted constant sheaf is perverse; this implies that the reduced cohomology
of the links of Whitney strata are concentrated in middle degree. As we discussed in Example 3.12, this allows us to
algebraically calculate the Betti numbers of the links. In the case of a general space, the obstruction to algebraically
calculating Milnor numbers is that there is no general algebraic manner for calculating the Betti numbers of the links
of strata.

\vskip .1in

Finally, we wish to say a few words about future directions for our work. In [{\bf Ma2}], [{\bf Ma3}], and [{\bf Ma4}], 
we developed the L\^e cycles and L\^e numbers of an affine hypersurface singularity. These L\^e numbers appear to be the
``correct'' generalization of the Milnor number to the case of arbitrary, non-isolated, affine hypersurface
singularities. Now, in this paper, we have generalized the Milnor number to the case of isolated hypersurface
singularities on an arbitrary analytic space. By combining these two approaches, we can obtain a super generalization of
the Milnor number -- one that works for arbitrary analytic functions on arbitrary analytic spaces. Moreover, using this
generalization, we can prove a super generalization of the result of L\^e and Saito in [{\bf L-S}]. Of course, as we
discussed above, the problem of actually calculating these generalized Milnor-L\^e numbers is precisely the problem of
calculating the Betti numbers of the complex links of Whitney strata.

\enddocument